\documentclass{article}
\usepackage{latexsym, amssymb,amsmath, amsbsy, amsopn, amstext, amsthm}
\usepackage{graphicx}
\usepackage{hyperref}
\usepackage{float}
\usepackage[textwidth=16cm,textheight=21.5cm]{geometry}

\usepackage{tikz}
\DeclareMathOperator{\rank}{rank}
 \DeclareMathOperator{\Col}{Col}
\DeclareMathOperator{\argmax}{argmax}
 \DeclareMathOperator{\lcm}{lcm}

\def\ra{\rightarrow}

\def\a{\alpha}
\def\b{\beta}
\def\d{\delta}

\def\D{\Delta}

\def\J{J}

\newcommand{\R}{{\mathbb R}}

\def\dsum{\mathop{\sum}\limits}

\newtheorem{thm}{Theorem}[section]

\newtheorem{cor}[thm]{Corollary}
\newtheorem{dfn}[thm]{Definition}
\newtheorem{prp}[thm]{Proposition}
\newtheorem{exa}[thm]{Example}
\newtheorem{rem}[thm]{Remark}

\begin{document}

\begin{center}

{\Large\bf Matrix Expression of Bayesian Game
\footnote{ Supported partly by NNSF 62073315 of China, and China Postdoctoral Science Foundation 2020TQ0184.} }
\end{center}

\vskip 2mm

\begin{center}
{\Large  Daizhan Cheng\dag, Changxi Li\ddag}

\vskip 2mm

\dag Institute of  Systems Science,
Chinese Academy of Sciences, Beijing 100190, P.R.China\\
E-mail: dcheng@iss.ac.cn\\

\ddag Shandong University\\

\vskip 2mm

\end{center}

\vskip \baselineskip

\underline{\bf Abstract}: A matrix-based framework for Bayesian games is presented, using semi-tensor product of matrices. Static Bayesian games are considered first. Matrix expression of Bayesian games is proposed. Three kinds of conversions, which convert Bayesian games to complete information games are investigated, certain properties are obtained, including two kinds of Bayesian-Nash equilibriums. Finally, dynamic Bayesian games are considered. Markoven dynamic equations are obtained for some strategy updating rules.

\vskip 5mm


\vskip 5mm

\section{Introduction}

The idea of game theory has a very long history. But it is a common believe that the symbol of modern game theory is the book of von Neumann and Morgenstern \cite{von44}. Modern game theory and modern control theory are like a twin. Not only because they were born in almost same time period, but also they have similar purpose: to ``manipulate"  objects according to players' purpose. Nowadays, more people try to use techniques in game theory to control problems, which are sometimes called game-based control, or game theoretic control \cite{aur14,bas98,fri12,tem10,yaz13,zha19}.

Bayesian game is also called the incomplete information game. Comparing with complete information game, which is also called normal noncooperative game,  Bayesian has more practical applications, because games with incomplete information, or uncertainty, are widely existing in real world, for instance \cite{gib92} gives some interesting examples such as Cournot competition under asymmetric information, first-price, sealed-bid auction, etc. Its various  applications are investigated \cite{heu96,giu09,eks15}.

Recently, semi-tensor product (STP) of matrices has been used to study some problems in game theory. For instance, there are STP-based investigations for the evolution of dynamic games \cite{gop11,che15}, potential games \cite{che14,liu16}, vector space structure of finite games \cite{che16,hao18}, and Boolean games \cite{che18}, just to mention a few.

Potential game was firstly proposed by \cite{ros73}, and then systematically developed in \cite{mon96}. It is of particular importance in game theoretic control \cite{can11}.  Various problems about Bayesian potential games have also been investigated \cite{ein20,eks16,eks17,fac97}.

The purpose of this paper is to provide a framework for Bayesian game. Roughly speaking, the major technique to deal with Bayesian game is to convert it into a normal game, which is also called complete information game. This method was firstly proposed by Harsanyi \cite{har67,har68,har68b}. It is recently called the Harsanyi transformation.
Another commonly used transformation is Selten transformation \cite{zam12}. Using STP, a very recent work proposed a new transformation, called Ex-ante agent transformation of Bayesian games, which keeps potential property unchanged \cite{wupr}.

The paper is organized as follows: Section 2 consists of some necessary preliminaries, including STP, normal noncooperative game, and potential game. Section 3 considers the matrix expression of static Bayesian game. In Section 4, two types fo Bayesian-Nash equilibriums are proposed. In Section 5, three kinds of conversions from Bayesian game to normal game. Section 6 considers two kinds of Bayesian potential games. Their verifications and some properties are presented. In Section 7, the dynamic Bayasian games are investigated. Particular interest has been put on some stratify updating rules, which lead to Markovin dynamic equations. Section 8 is a brief conclusion.

\section{Preliminaries}

\subsection{STP of matrices}

This subsection briefly reviews STP of matrices.

\begin{dfn} \label{d2.1.1} \cite{che11, che12a}:
Let $M\in {\cal M}_{m\times n}$, $N\in {\cal M}_{p\times q}$, and $t=\lcm\{n,p\}$ be the least common multiple of $n$ and $p$.
The semi-tensor product (STP) of $M$ and $N$, denoted by $M\ltimes N$, is defined as
\begin{align}\label{2.1.1}
\left(M\otimes I_{t/n}\right)\left(N\otimes I_{t/p}\right)\in {\cal M}_{mt/n\times qt/p},
\end{align}
where $\otimes$ is the Kronecker product.
\end{dfn}

Note that when $n=p$, $M\ltimes N=MN$. That is, the semi-tensor product is a generalization of conventional matrix product. Moreover, it keeps all the properties of conventional matrix product available \cite{che12a}. Hence we can omit the symbol $\ltimes$. Throughout this paper the matrix product is assumed to be STP, and the symbol $\ltimes$ is mostly omitted.

The following are some basic properties:

\begin{prp}\label{p2.1.2}
\begin{enumerate}
\item (Associative Law)
\begin{align}\label{2.1.2}
(F\ltimes G)\ltimes H = F\ltimes (G\ltimes H).
\end{align}
\item (Distributive Law)
\begin{align}\label{2.1.3}
\begin{cases}
F\ltimes (aG\pm bH)=aF\ltimes G\pm bF\ltimes H,\\
(aF\pm bG)\ltimes H=a F\ltimes H \pm bG\ltimes H,\quad a,b\in \R.
\end{cases}
\end{align}
\end{enumerate}
\end{prp}

\begin{prp}\label{p2.1.3}
\begin{enumerate}
\item Let $X\in \R^m$, $Y\in \R^n$ be two columns. Then
\begin{align}\label{2.1.4}
X\ltimes Y=X\otimes Y.
\end{align}
\item Let $\omega \in \R^m$, $\sigma\in \R^n$ be two rows. Then
\begin{align}\label{2.1.5}
\omega\ltimes \sigma =\sigma \otimes \omega.
\end{align}
\end{enumerate}
\end{prp}

The following proposition shows that the STP satisfies the block-multiplication law.

\begin{prp}\label{p2.1.4}
Let
$$
A=\begin{bmatrix}
A^{11}&A^{12}&\cdots&A^{1q}\\
A^{21}&A^{22}&\cdots&A^{2q}\\
\vdots&~&~&~\\
A^{p1}&A^{p2}&\cdots&A^{pq}\\
\end{bmatrix}
,\quad
B=\begin{bmatrix}
B^{11}&B^{12}&\cdots&B^{1r}\\
B^{21}&B^{22}&\cdots&B^{2r}\\
\vdots&~&~&~\\
B^{q1}&B^{q2}&\cdots&B^{qr}\\
\end{bmatrix},
$$
where
$$
\begin{array}{l}
n_c(A^{ij})=n_c(A^{ik}),\quad 1\leq j,k\leq q,\;\forall i,\\
n_r(B^{ik})=n_r(B^{jk}),\quad 1\leq i,j\leq q,\;\forall k.
\end{array}
$$
($n_c(A)$ ($n_r(A)$) is the column (row) of $A$.
Then
$$
A\ltimes B=\left(C^{i,j}\right),
$$
where
$$
C^{i,j}=\dsum_{k=1}^q A^{i,k}\ltimes B^{k,j},\quad i=1,2,\cdots,p;\;j=1,2,\cdots, r.
$$
\end{prp}

About the transpose, we have
\begin{prp}\label{p2.1.5}
\begin{align}
\label{2.1.6} (A\ltimes B)^\mathrm{T}=B^\mathrm{T}\ltimes A^\mathrm{T}.
\end{align}
\end{prp}

About the inverse, we have

\begin{prp}\label{p2.1.6}
Assume $A$ and $B$ are invertible, then
\begin{align} \label{2.1.7}
(A\ltimes B)^{-1}=B^{-1}\ltimes A^{-1}.
\end{align}
\end{prp}

The following property is for STP only.

\begin{prp}\label{p2.1.7} Let $X\in \R^m$ be a column and $M$ a matrix. Then
\begin{align}\label{2.1.8}
X\ltimes M=\left(I_m\otimes M\right)X.
\end{align}
\end{prp}

Let $m\geq 1$, $n\geq 1$, and $1\leq k\leq mn$. Define
\begin{align}\label{2.1.8}
j_k:=k(mod~n);\quad i_k=(k-j)/n+1.
\end{align}

\begin{dfn}\label{d2.1.8} \cite{che11} A matrix $W_{[m,n]}\in {\cal M}_{mn\times mn}$, defined by
\begin{align}\label{2.1.9}
\Col_k\left(W_{[m,n]}\right):=\d_{mn}^{(j_k-1)m+i},\quad k=1,\cdots,mn,
\end{align}
is called the $(m,n)$-th dimensional swap matrix.
\end{dfn}

The basic function of the swap matrix is to ``swap" two vectors. That is,

\begin{prp}\label{p2.1.10} Let $X\in \R^m$ and $Y\in \R^n$ be two columns. Then
\begin{align}\label{2.1.11}
W_{[m,n]}\ltimes X\ltimes Y=Y\ltimes X.
\end{align}
\end{prp}

Define a matrix $PR_n\in {\cal L}_{n^2\times n}$, called the order reducing matrix, as
\begin{align}\label{2.1.16}
\Col_i(PR_n)=\d_{n^2}^{(i-1)n+i},\quad i=1,\cdots,n.
\end{align}

\begin{dfn}
Let $A\in \mathcal{M}_{p\times n}$ and $B\in \mathcal{M}_{q\times n}$. Then the Khatri-Rao Product of $A$ and $B$ is
  \begin{align}
  \begin{array}{ccl}
  A*B=[\Col_1(A)\ltimes \Col_1(B),\cdots,\Col_n(A)\ltimes \Col_n(B)]
  \in \mathcal{M}_{pq\times n}.
  \end{array}
  \end{align}
\end{dfn}

\subsection{Normal Non-cooperative Game}

\begin{dfn}\label{d2.2.1}  A (finite non-cooperative) normal  game is described as a triple $G=(N,S,C)$, where
\begin{itemize}
\item[(i)] $N=\{1,2,\cdots,n\}$ is the set of players;
\item[(ii)] $S=\prod_{i=1}^nS_i$ is called the profile, with
\begin{align}\label{2.2.1}
S_i=\{s^i_1,s^i_2,\cdots,s^i_{k_i}\},\quad i=1,\cdots,n,
\end{align}
the set of strategies of player $i$, and  $S_{-i}:=\prod_{j\neq i}^nS_j$;
\item[(iii)] $C=(c_1,\cdots,c_n)\in \R^n$ with $c_i:S\ra \R$ defined as
\begin{align}\label{2.2.2}
c_i:=c_i(s_1,\cdots,s_n), \quad s_j\in S_j,\; j=1,\cdots,n,\quad i=1,\cdots,n,
\end{align}
is called the payoff function of player $i$.
\end{itemize}
\end{dfn}

Denote the structure vector of $c_i$ by $V_i^c$, $i=1,\cdots,n$. Assume $|N|=n$ and $k_i$, $i=1,\cdots,n$ are fixed, we denote the set of such games by
${\cal G}_{[n;k_1,\cdots,k_n]}$. Note that if $G\in {\cal G}_{[n;k_1,\cdots,k_n]}$, then $G$ is uniquely determined by $\{V_i^c\;\big|\; i=1,\cdots,n\}$. Hence,
\begin{align}\label{2.2.3}
{\cal G}_{[n;k_1,\cdots,k_n]}\simeq \R^{n\kappa},
\end{align}
where $\kappa=\prod_{i=1}^nk_i$.

\begin{dfn}\label{d2.2.2} $G\in {\cal G}_{[n;k_1,\cdots,k_n]}$ is expressed as $V_G=(V_1^c,\cdots,V_n^c)\in \R^{n\kappa}$, where $V_i^c$ is the structure vector of $c_i$ of $G$, $i=1,\cdots,n$.
\end{dfn}

\begin{rem}\label{r2.2.3}  The decision-making (or strategy-selection of each player $j$) is based on the difference of $c_j$ over two different strategy profiles
$$
c_j(s)-c_j(s'), \quad j=1,\cdots,n;\;s,s'\in S.
$$
For instance, the following concepts are depending on such differences:
\begin{enumerate}

\item[(i)] A strategy profile $s^*=(s_1^*,\cdots,s_n^*)$ is a Nash equilibrium if
\begin{align}\label{2.2.4}
c_i(s_i^*,s_{-i}^*)\geq c_i(s_i,s_{-i}^*),\quad\forall s_i\in S_i,~s_{-i}^*\in S_{-i}, i=1,\cdots,n.
\end{align}
It is easy to verify that Nash equilibrium depends on $\{c_j(s)-c_j(s')\;\big|\;j=1,\cdots,n\}$ only.

\item[(ii)] In the definition of potential games \cite{mon96} it is assumed that there is a pseudo-mixed-logical function $P:\prod_{i=1}^n{\cal D}_{k_i}\ra \R$ such that
$$
c_i(x_i,s_{-i})-c_i(y_i,s_{-i})=P(x_i,s_{-i})-P(y_i,s_{-i}),\quad i=1,\cdots,n.
$$
It is clear that whether a game is potential depends on $\{c_j(s)-c_j(s')\;\big|\;j=1,\cdots,n\}$ only.

\item[(iii)] Consider an evolutionary game, two of the most commonly used strategy updating rules (SURs) are the better and best response dynamics \cite{oc13}.
\begin{itemize}
\item Better response SUR: updating its strategy to a strategy in $\{s_i\;\big|\;c_i(s_i,s_{-i}(t))>c_i(s(t))\}$, chosen uniformly at random;
\item Best response SUR: updating its strategy to a strategy in $\argmax_{s_i}c_i(s_i,s_{-i}(t))$, chosen uniformly at random;
\end{itemize}
One sees easily that they depend on $\{c_j(s)-c_j(s')\;\big|\;j=1,\cdots,n\}$ only.
\end{enumerate}
\end{rem}

The following proposition is obvious.

\begin{prp}\label{p2.2.4} Two games $G$ and $G'$ are equivalent (for decision making) if there exists a row vector $v\in \R^{\kappa}$ such that
\begin{align}\label{2.2.5}
V^c_i+v=V^{c'}_i,\quad i=1,\cdots,n.
\end{align}
\end{prp}

From Proposition \ref{p2.2.4} one sees easy that a game (precisely, an equivalent class of games) can be expressed as
\begin{align}\label{2.2.6}
W_G:=[V^c_2 - V^c_1,\; V^c_3 - V^c_1,\; \cdots,\; V^c_n - V^c_1].
\end{align}
We call (\ref{2.2.6}) the equivalent vector form of $G$.

Using equivalent vector form of $G$, one sees that $G_{[n;k_1,\cdots,k_n]}$ is an $(n-1)\kappa$ dimensional vector space. That is,
\begin{align}\label{2.2.7}
G_{[n;k_1,\cdots,k_n]}\sim \R^{(n-1)\kappa}.
\end{align}

\subsection{Potential Game}

\begin{dfn} \cite{mon96} \label{d2.3.1} Consider a finite game $G\in {\cal G}_{[n;k_1,\cdots,k_n]}$.
If there is a function $P:S\ra \R$, such that for every $i\in N$ and for every $s_{-i}\in S_{-i}$
\begin{align}\label{2.3.1}
c_i(x, s_{-i})-c_i(y,s_{-i})=P(x, s_{-i})-P(y,s_{-i}), \quad i=1,\cdots,n,
\end{align}
then $G$ is called an exact potential game.
\end{dfn}

The followings are some fundamental properties of potential.

\begin{thm}\label{t2.3.2} \cite{mon96} If $G$ is a potential game, then the potential function $P$ is unique up to a constant number. Precisely, if $P_1$ and $P_2$ are two potential functions, then $P_1-P_2=c_0\in \R$.
\end{thm}

\begin{thm}\label{t2.3.3} \cite{mon96} Let $P$ be a potential function for $G$. Then $s\in S$ is an equilibrium point of $G$, if and only if
\begin{align}\label{2.3.2}
P(s)\geq P(s_{-i},x),\quad \forall x\in S_i,\; i=1,\cdots,n.
\end{align}
Particularly, if $P$ admits a maximal value in $S$, then $G$ has at least one pure Nash equilibrium.
\end{thm}

\begin{cor}\label{c2.3.4} Every finite potential game possesses a pure Nash equilibrium.
\end{cor}

The following argument is a generalization of potential equation \cite{che14}.

Assume $G\in {\cal G}_{[n;k_1,\cdots,k_n]}$ is potential, define
\begin{align}\label{2.3.3}
d_i(s_1,\cdots,s_n):=c_i(s_1,\cdots,s_n)-P(s_1,\cdots,s_n),\quad i=1,\cdots,n.
\end{align}
Then it is easy to verify that $d_i$ is independent of $s_i$, denoted by $d_i=d_i(s_1,\cdots,\hat{s}_i,\cdots,s_n)$.
In vector form (\ref{2.3.3}) becomes
\begin{align}\label{2.3.4}
V^d_i\ltimes_{j\neq i}s_j=V^c_is-V^Ps,\quad i=1,\cdots,n.
\end{align}
where $s=\ltimes_{j=1}^ns_j$.

Define
$$
E_i:=I_{\a_i}\otimes {\bf 1}^T_{k_i}\otimes I_{\b_i},\quad i=1,\cdots,n,
$$
where
$$
\begin{array}{ll}
\a_1=1,&\a_i=\prod_{j=1}^{i-1}k_j,\; j\geq 2\\
\b_n=1,&\b_i=\prod_{j=i+1}^{n}k_j,\; j\leq n-1.
\end{array}
$$
Then (\ref{2.3.4}) can be expressed as
\begin{align}\label{2.3.5}
V^d_i E_i=V^c_i-V^p,\quad i=1,\cdots,n.
\end{align}
Solving the first equation to get
$$
V^P=V^c_1-V^d_1E_1,
$$
and plugging it into the other equations of (\ref{2.3.5}) yield
\begin{align}\label{2.3.6}
V^d_i E_i=V^d_1E_1=V^c_i-V^c_1,\quad i=2,3,\cdots,n.
\end{align}
Denote
\begin{align}\label{2.3.7}
\xi_i=(V^d_i)^T,\quad i=1,\cdots,n;\quad b:=W_G^T.
\end{align}
Then (\ref{2.3.6}) can be expressed as
\begin{align}\label{2.3.8}
\Psi\xi=b,
\end{align}
where $\xi=(\xi_1^T,\cdots,\xi_n^T)^T$ and
\begin{align}\label{2.3.9}
\Psi=\begin{bmatrix}
-E_1^T&E_2^T&0&\cdots&0\\
-E_1^T&0&E_3^T&\cdots&0\\
\vdots&~&~&~&~\\
-E_1^T&0&0&\cdots&E_n^T.
\end{bmatrix}
\end{align}

Summarizing the above argument, we have

\begin{thm}\label{t2.3.5} \cite{che16} Let $G\in {\cal G}_{[n,k_1,\cdots,k_n]}$. $G$ is a potential game, if and only if (\ref{2.3.8}) has solution.
Moreover, when a solution exists, then
\begin{align}\label{2.3.10}
V^P=V^c_1-\xi_1^TE_1.
\end{align}
\end{thm}

\section{Static Bayesian Game}

\begin{dfn}\label{d3.1}  \cite{gib92} A (finite) static Bayesian game (BG) $G=(N,T,A,c, Pr)$
consists of
\begin{itemize}
\item[(i)] Player set:
\begin{align}\label{3.1}
N=\{1,2,\cdots,n\}.
\end{align}

\item[(ii)] Type set:
\begin{align}\label{3.2}
T=T_1\times T_2\times\cdots\times T_n,
\end{align}
where $T_i$ is the type set of player $i$ with
\begin{align}\label{3.3}
T_i=\{t_i^1,t_i^2,\cdots,t_i^{\tau_i}\},\quad i=1,2,\cdots,n.
\end{align}

\item[(iii)] Profile and actions:
\begin{align}\label{3.4}
A=A_1\times A_2\times \cdots\times A_n,
\end{align}
where $A_i$ is the action set of player $i$ with
\begin{align}\label{3.5}
A_i=\left\{a_i^1,a_i^2,\cdots,a_i^{r_i}\right\},\quad i=1,2,\cdots,n,
\end{align}
and $A_i(t_i^j)$ is the admissible actions under type $t_i^j$ with
\begin{align}\label{3.6}
A_i(t_i^j)\subset A_i,\quad j=1,2,\cdots,\tau_i,\; i=1,\cdots,n.
\end{align}
\item[(iv)] Payoff functions:
\begin{align}\label{3.7}
c_i:A_1(t_1^{j_1})\times A_2(t_2^{j_2}) \times \cdots\times A_n(t_n^{j_n})\ra \R,\quad i=1,2,\cdots,n.
\end{align}

\item[(v)]  Beliefs:
\begin{align}\label{3.8}
p_i\left(t_{-i}\;|\; t_i\right)=\dfrac{Pr(t_1,t_2,\cdots,t_n)}{Pr(t_i)}=\dfrac{Pr(t_i,t_{-i})}{\dsum_{t'_{-i}}Pr(t_i,t'_{-i})},
\end{align}
where $t_{-i}\in T_{-i}:=\prod_{j\neq i}T_j$. $Pr(t_1,t_2,\cdots,t_n)$, called the common prior, is a common knowledge, and $t_i$ is a private knowledge for player $i$.

\end{itemize}
\end{dfn}

Identify  $i\sim \d_l^i\in D_l$, then $\d_l^i$ is called the vector representation of integer $i.$ By virtue of vector representation of actions and types, then the action set $A_i$ and type set $T_i$ can be expressed as
$$
A_i=\left\{\d_{r_i}^j\;|\;1\leq j\leq r_i\right\};
T_i=\left\{\d_{\tau_i}^j\;|\;1\leq j\leq \tau_i\right\}, \quad i=1,2,\cdots,n.
$$
Denote by $r=\prod_{i=1}^nr_i$, $\tau=\prod_{i=1}^n\tau_i$, $R_{-\infty}=R\cup \{-\infty\}$, and
\begin{align}\label{3.9}
\begin{array}{l}
\bar{c}_i(a_1,\cdots,a_n;t_1,\cdots,t_n)
=\begin{cases}
c_i(a_1,\cdots,a_n;t_1,\cdots,t_n),\quad a_j\in A_j(t_j),\\
-\infty,\quad \mbox{otherwise},\quad 1\leq i\leq n.
\end{cases}
\end{array}
\end{align}
Then we have
\begin{align}\label{3.10}
\bar{c}_i:\D_{r\tau}\ra \R_{-\infty},\quad 1\leq i\leq n.
\end{align}

Using (\ref{3.9}), it is obvious that for each $i$ there exists a unique row vector $V_i\in \R_{-\infty}^{r\tau}$ such that
\begin{align}\label{3.11}
\bar{c}_i(a_1,\cdots,a_n;t_1,\cdots,t_n)=V_i^{\bar{c}}ta,\quad 1\leq i\leq n,
\end{align}
where $a=\ltimes_{i=1}^na_i$, $t=\ltimes_{i=1}^nt_i$. Let $V_G:=[V_1^{\bar{c}},\cdots,V_n^{\bar{c}}]$, which is called the structure vector of Bayesian game $G$.

\begin{rem}\label{r3.2}
\begin{itemize}
\item[(i)] Here we assume the purpose of each player is to maximize his payment. So using $\bar{c}_i$ to optimize payment is equivalent to $c_i$. If the purpose of each player is to minimize his cost, then in (\ref{3.9}) $+\infty$ is used to replace $-\infty$.
\item[(ii)] In fact, (\ref{3.11}) provides a natural vector space structure for $c_i(a,t)\in \R^{r\tau}$, where $a=a_1a_2\cdots a_n\in \D_r$, $t=t_1t_2\cdots t_n\in \D_\tau$.
\end{itemize}
\end{rem}

\begin{exa}\label{e3.3} Consider a finite BG
\begin{align}\label{3.12}
G=(N,T,A,c,p),
\end{align}
where
$$
\begin{array}{ccl}
N&=&\{1,2\};\\
T&=&\{T_1,T_2\},\;\mbox{with}\\
T_1&=&\{t_1^1,t_1^2\},\\
T_2&=&\{t_2^1,t_2^2\};\\
A&=&\{A_1,A_2\},\;\mbox{with}\\
A_1&=&\{a_1^1,a_1^2,a_1^3\},\\
A_2&=&\{a_2^1,a_2^2,a_2^3\},\\
A_1(t_1^1)&=&\{a_1^1,a_1^2\},\\
A_1(t_1^2)&=&\{a_1^2,a_1^3\},\\
A_2(t_2^1)&=&\{a_2^1,a_2^2,a_2^3\},\\
A_2(t_2^2)&=&\{a_2^1,a_2^3\}.\\
\end{array}
$$

\begin{table}[!htb]
\centering
\caption{Payment Function for $t_1^1-t_2^1$}\label{tb.1.1}
\vskip .25\baselineskip
\begin{tabular}{|c||c|c|c|}
\hline
$c_1\backslash c_2$&$a_2^1$&$a_2^2$&$a_2^3$\\
\hline
\hline
$a_1^1$&$2,~3$&$1,~4$&$1,~-2$\\
\hline
$a_1^2$&$1,~-2$&$2,~1$&$0,-3$\\
\hline
\end{tabular}
\end{table}

\begin{table}[!htb]
\centering
\caption{Payment Function for $t_1^1-t_2^2$}\label{tb.1.2}
\vskip .25\baselineskip
\begin{tabular}{|c||c|c|}
\hline
$c_1\backslash c_2$&$a_2^1$&$a_2^3$\\
\hline
\hline
$a_1^1$&$-1,~2$&$1,~3$\\
\hline
$a_1^2$&$1,~-2$&$-2,~0$\\
\hline
\end{tabular}
\end{table}

\begin{table}[!htb]
\centering
\caption{Payment Function for $t_1^2-t_2^1$}\label{tb.1.3}
\vskip .25\baselineskip
\begin{tabular}{|c||c|c|c|}
\hline
$c_1\backslash c_2$&$a_2^1$&$a_2^2$&$a_2^3$\\
\hline
\hline
$a_1^2$&$3,~5$&$2,~4$&$2,~0$\\
\hline
$a_1^3$&$2,~-2$&$-2,~4$&$3,~3$\\
\hline
\end{tabular}
\end{table}

\begin{table}[!htb]
\centering
\caption{Payment Function for $t_1^2-t_2^2$}\label{tb.1.4}
\vskip .25\baselineskip
\begin{tabular}{|c||c|c|}
\hline
$c_1\backslash c_2$&$a_2^1$&$a_2^3$\\
\hline
\hline
$a_1^2$&$2,~1$&$-1,~-3$\\
\hline
$a_1^3$&$2,~2$&$-1,~-2$\\
\hline
\end{tabular}
\end{table}

\begin{table}[!htb]
\centering
\caption{Distribution}\label{tb.1.5}
\vskip .25\baselineskip
\begin{tabular}{|c||c|c|}
\hline
$t_1\backslash t_2$&$t_2^1$&$t_2^2$\\
\hline
\hline
$t_1^1$&$0.3$&$0.2$\\
\hline
$t_1^2$&$0.1$&$0.4$\\
\hline
\end{tabular}
\end{table}

Then
$$
\begin{array}{l}
P_r(t^1_2\;|\;t_1^1)=\frac{P_r(t_1^1\cup t_2^1)}{P_r(t^1_1)}
=\frac{P_r(t_1^1\cup t_2^1)}{P_r(t^1_1\cap t^1_2)+P_r(t^1_1\cap t^2_2)}=0.75\\
P_r(t^2_2\;|\;t_1^1)=0.25;\quad
P_r(t^1_2\;|\;t_1^2)=1/3;\quad
P_r(t^2_2\;|\;t_1^2)=2/3\\
P_r(t^1_1\;|\;t_2^1)=0.6;\quad
P_r(t^2_1\;|\;t_2^1)=0.4;\quad
P_r(t^1_1\;|\;t_2^2)=0.2\\
P_r(t^2_1\;|\;t_2^2)=0.8\\
\end{array}
$$
Hence the beliefs are:
\begin{align}\label{3.13}
\begin{array}{l}
p_{t_1^1}=(0.75,~0.25)^T,\quad p_{t_1^2}=(1/3,~2/3)^T\\
p_{t_2^1}=(0.6,~0.4)^T,\quad p_{t_2^2}=(0.2,~0.8)^T.
\end{array}
\end{align}

\begin{table}[!htb]
\centering
\caption{Payoffs}\label{tb.1.6}
\vskip .25\baselineskip
\begin{tabular}{|c||c|c|c|c|}
\hline
$c\backslash (t,a)$&$(\d_4^1,\d_9^1)$&$(\d_4^1,\d_9^2)$&$(\d_4^1,\d_9^3)$&$(\d_4^1,\d_9^4)$\\
\hline
\hline
$\bar{c}_1$&$2$&$1$&$1$&$1$\\
\hline
$\bar{c}_2$&$3$&$4$&$-2$&$-2$\\
\hline
\hline
$c\backslash (t,a)$&$(\d_4^1,\d_9^5)$&$(\d_4^1,\d_9^6)$&$(\d_4^1,\d_9^7)$&$(\d_4^1,\d_9^8)$\\
\hline
\hline
$\bar{c}_1$&$2$&$0$&$-\infty$&$-\infty$\\
\hline
$\bar{c}_2$&$1$&$-3$&$-\infty$&$-\infty$\\
\hline
\hline
$c\backslash (t,a)$&$(\d_4^1,\d_9^9)$&$(\d_4^2,\d_9^1)$&$(\d_4^2,\d_9^2)$&$(\d_4^2,\d_9^3)$\\
\hline
\hline
$\bar{c}_1$&$-\infty$&$-1$&$-\infty$&$1$\\
\hline
$\bar{c}_2$&$-\infty$&$2$&$-\infty$&$3$\\
\hline
\end{tabular}
\end{table}

\begin{table}[!htb]
\centering
\caption{Payoffs(Cont'd)}
\vskip .25\baselineskip
\begin{tabular}{|c||c|c|c|c|}
\hline
$c\backslash (t,a)$&$(\d_4^2,\d_9^4)$&$(\d_4^2,\d_9^5)$&$(\d_4^2,\d_9^6)$&$(\d_4^2,\d_9^7)$\\
\hline
\hline
$\bar{c}_1$&$1$&$-\infty$&$-2$&$-\infty$\\
\hline
$\bar{c}_2$&$-2$&$-\infty$&$0$&$-\infty$\\
\hline
\hline
$c\backslash (t,a)$&$(\d_4^2,\d_9^8)$&$(\d_4^2,\d_9^9)$&$(\d_4^3,\d_9^1)$&$(\d_4^3,\d_9^2)$\\
\hline
\hline
$\bar{c}_1$&$-\infty$&$-\infty$&$-\infty$&$-\infty$\\
\hline
$\bar{c}_2$&$-\infty$&$-\infty$&$-\infty$&$-\infty$\\
\hline
\hline
$c\backslash (t,a)$&$(\d_4^4,\d_9^3)$&$(\d_4^3,\d_9^4)$&$(\d_4^3,\d_9^5)$&$(\d_4^3,\d_9^6)$\\
\hline
\hline
$\bar{c}_1$&$-\infty$&$3$&$2$&$2$\\
\hline
$\bar{c}_2$&$-\infty$&$5$&$4$&$0$\\
\hline
\end{tabular}
\end{table}

\begin{table}[!htb]
\centering
\caption{Payoffs(Cont'd)}
\vskip .25\baselineskip
\begin{tabular}{|c||c|c|c|c|}
\hline
$c\backslash (t,a)$&$(\d_4^3,\d_9^7)$&$(\d_4^3,\d_9^8)$&$(\d_4^3,\d_9^9)$&$(\d_4^4,\d_9^1)$\\
\hline
\hline
$\bar{c}_1$&$2$&$-2$&$3$&$-\infty$\\
\hline
$\bar{c}_2$&$-2$&$4$&$3$&$-\infty$\\
\hline
\hline
$c\backslash (t,a)$&$(\d_4^4,\d_9^2)$&$(\d_4^4,\d_9^3)$&$(\d_4^4,\d_9^4)$&$(\d_4^4,\d_9^5)$\\
\hline
\hline
$\bar{c}_1$&$-\infty$&$-\infty$&$2$&$-\infty$\\
\hline
$\bar{c}_2$&$-\infty$&$-\infty$&$-1$&$-\infty$\\
\hline
\hline
$c\backslash (t,a)$&$(\d_4^4,\d_9^6)$&$(\d_4^4,\d_9^7)$&$(\d_4^4,\d_9^8)$&$(\d_4^4,\d_9^9)$\\
\hline
\hline
$\bar{c}_1$&$1$&$2$&$-\infty$&$-1$\\
\hline
$\bar{c}_2$&$-3$&$2$&$-\infty$&$-2$\\
\hline
\end{tabular}
\end{table}

The cost functions are described as follows:

$$
\begin{cases}
\bar{c}_1(a,t)=V^{\bar{c}}_1ta,\\
 \bar{c}_2(a,t)=V^{\bar{c}}_2ta,
\end{cases}
$$
where
$$
\begin{array}{ccl}
V^{\bar{c}}_1&=&[2,1,1,1,2,0,-\infty,-\infty,-\infty,-1,-\infty,1,1,\\
~&~&-\infty,-2,-\infty,-\infty,-\infty,
-\infty,-\infty,-\infty,3,2,2,\\
~&~&2,-2,3,-\infty,-\infty,-\infty,2,-\infty,-1,2,-\infty,-1]\\
V^{\bar{c}}_2&=&[3,4,-2,-2,1,-3,-\infty,-\infty,-\infty,2,-\infty,3,-2,\\
~&~&-\infty,0,-\infty,-\infty,-\infty,
-\infty,-\infty,-\infty,5,4,0,\\
~&~&-2,4,3,-\infty,-\infty,-\infty,-1,-\infty,3,2,-\infty,-2]\\
\end{array}
$$
\end{exa}

\section{Bayesian-Nash Equilibrium}

This section proposes two kinds of types first.

\begin{dfn}\label{d4.1} There are two kinds of types:
\begin{itemize}
\item[(i)] Types of Nature (TN):

The types are determined by pre-assigned distribution $Pr(t_1,t_2,\cdots,t_n)$, which is a common knowledge.
Each player $i$ knows the type $t_i$ assigned to him. The type is assigned  by Nature.

\item[(ii)] Types of Human (TH):

Player $i$ has the right to choose $t_i$ for entering the game. Then $t_i$ becomes part of strategy for player $i$, $i=1,2,\cdots, n$.

\end{itemize}
\end{dfn}

The expected payoff functions under different  kinds of types are defined as follows.
\begin{prp}\label{p4.2}

\begin{itemize}
\item[(i)] TN:

The expected value of player $i$, using actions $a=(a_1,a_2,\cdots,a_n)$ is

\begin{align}\label{4.1}
\begin{array}{ccl}
E^N_i(a)&:=&E_i(a(T))\\
&=&\dsum_{t\in T}Pr(t)\bar{c}_i(a,t),\\
&=&\dsum_{t\in T}Pr(t)V^{\bar{c}}_i(a,t),\quad i=1,2,\cdots,n.
\end{array}
\end{align}

\item[(ii)] TH:

The expected value of player $i$, using actions $a=(a_1,a_2,\cdots,a_n)$ and type $t_i^j$ is
\begin{align}\label{4.2}
\begin{array}{ccl}
E^H_i(a,t_i^j)&:=&E_i(a(T)|t_i=t_i^j)\\
~&=&\dsum_{t_{-i}\in T_{-i}}p_{i}(t_{-i}|t_i^j)\bar{c}_i(a,t^j_i,t_{-i}),\\
~&=&\dsum_{t_{-i}\in T_{-i}}p_{i}(t_{-i}|t_i^j)V^{\bar{c}}_ita,\quad i=1,2,\cdots,n.
\end{array}
\end{align}
\end{itemize}
\end{prp}
\noindent{\it Proof:}
It is easy to verify the results by virtue of (\ref{3.11})
\hfill $\Box$

\begin{dfn}\label{d4.3} Consider a  (finite) static Bayesian game $G=(N,T,A,c, Pr)$.

\begin{itemize}
\item[(i)] TN:

A profile $\left(a_1^*,a_2^*,\cdots,a_n^*\right)$ is said to be a pure Bayesian-Nash equilibrium, if for each $i$ and any $t_i^j\in T_i$ the following inequalities hold.
\begin{align}\label{4.3}
\begin{array}{l}
E^N_i(a^*(t)|t_i)\geq E^N_i\left(a_1^*(t_1),\cdots, a_{i-1}^*(t_{i-1}), a_i,a_{i+1}^*(t_{i+1}),\cdots, a_{n}^*(t_n)|t_i\right),
\quad \forall t_i\in T_i;\;i=1,2,\cdots,n.
\end{array}
\end{align}

\item[(ii)] TH:

A profile $\left(a_1^*(t_1^*),a_2^*(t_2^*),\cdots,a_n^*(t_n^*)\right)$  is said to be a pure Bayesian-Nash equilibrium, if for each $i$ and any $t\in T$ the following inequalities hold.
\begin{align}\label{2.4}
\begin{array}{l}
E^H_i\left(a_1^*(t_1^*),a_2^*(t_2^*),\cdots,a_n^*(t_n^*)\right)\geq E^H_i\left(a_1^*(t_1^*),\cdots, a_{i-1}^*(t^*_{i-1}),\right.\\
~~\left.a_i(t_i), a_{i+1}^*(t_{i+1}^*),\cdots, a_{n}^*(t_n^*)\right), \forall t_i\in T_i,\; i=1,2,\cdots,n.
\end{array}
\end{align}
\end{itemize}
\end{dfn}

\begin{dfn}\label{d4.4} Consider a  (finite) static Bayesian game $G=(N,T,A,c, Pr)$.
 The belief vector of player $i$ for $t_i^j$  is defined as
\begin{align}\label{4.5}
\begin{array}{lllccc}
p_i(t_i^j):&=&Pr(t_{-i}\in T_{-i}|t_i=t_i^j)~~~~~~~~~\\
&=&\begin{bmatrix}
Pr\left(t_{-i}=(t_1^1,\cdots,t_{i-1}^1,t_{i+1}^1\cdots,t_n^1)|t_i=t_i^j\right)\\
Pr\left(t_{-i}=(t_1^1,\cdots,t_{i-1}^1,t_{i+1}^1\cdots,t_n^2)|t_i=t_i^j\right)\\
\cdots\\
Pr\left(t_{-i}=(t_1^{\tau_1},\cdots,t_{i-1}^{\tau_{i-1}},t_{i+1}^{\tau_{i+1}}\cdots,t_n^{\tau_n})|t_i=t_i^j\right)\\
\end{bmatrix},\\
&~~&~~~~~\\
~~~~~p_i&=&[p_i(t_i^1),p_i(t_i^2),\cdots,p_i(t_i^{\tau_i})],~~~~i=1,2,\cdots,n.
\end{array}
\end{align}
\end{dfn}

\begin{exa}\label{4.5} Consider a Bayesian game, which has distribution as in Table \ref{tb.4.7}.

\begin{table}[h]
\centering
\caption{Distribution}\label{tb.4.7}
\vskip .25\baselineskip
\begin{tabular}{|c||c|c|c|}
\hline
$t_1\backslash t_2$&$t_2^1$&$t_2^2$&$t_2^3$\\
\hline
\hline
$t_1^1$&$0.1$&$0.2$&$0.3$\\
\hline
$t_1^2$&$0.15$&$0.1$&$0.15$\\
\hline
\end{tabular}
\end{table}

$$
p_1({t_1^1})=\begin{bmatrix}
Pr(t_2^1|t_1^1)\\
Pr(t_2^2|t_1^1)\\
Pr(t_2^3|t_1^1)\\
\end{bmatrix}
=\begin{bmatrix}
\frac{Pr(t_2^1,t_1^1)}{Pr(t_2^1,t_1^1)Pr(t_2^2,t_1^1)Pr(t_2^3,t_1^1)}\\
~~\\
\frac{Pr(t_2^2,t_1^1)}{Pr(t_2^1,t_1^1)Pr(t_2^2,t_1^1)Pr(t_2^3,t_1^1)}\\
~~\\
\frac{Pr(t_2^3,t_1^1)}{Pr(t_2^1,t_1^1)Pr(t_2^2,t_1^1)Pr(t_2^3,t_1^1)}\\
\end{bmatrix}
=\begin{bmatrix}
1/6\\
1/3\\
1/2\\
\end{bmatrix}.
$$
Similarly,
$$
\begin{array}{ll}
p_1({t_1^2})=(3/8,1/4,3/8)^T,&p_2({t_2^1})=(0.4,0.6)^T,\\
p_2({t_2^2})=(2/3,1/3)^T,&p_2({t_2^3})=(2/3,1/3)^T.\\
\end{array}
$$
\end{exa}

\section{Conversions in Bayesian Games}
The key idea to deal with Bayesian game is to convert it into a normal game, which is also called complete information game. This method was firstly proposed by Harsanyi \cite{har67,har68,har68b}, which  is called Harsanyi conversion.
Another commonly used transformation is Selten conversion \cite{zam12}. Besides Harsanyi conversion and Selten conversion, we provide a new kind of conversion, called Action-Type conversion.

\begin{dfn}\label{d4.6} Consider a finite Bayesian game $G$ with prior probability distribution $Pr(t_1,t_2,\cdots,t_n)$.  Three conversions are defined as follows:
\begin{itemize}

\item[(i)] Harsanyi Conversion (H-Conversion):
Define
\begin{align}\label{4.6}
c_i^H(a):=Ec_i(a),\quad i=1,2,\cdots,n.
\end{align}

\item[(ii)] Selten Conversion (S-Conversion):  Player $i$ knows his type $t_i=\bar{t}_i$.
Define
\begin{align}\label{4.7}
c_i^S(a):=E(c_i(a)|t_i=\bar{t}_i),\quad i=1,2,\cdots,n.
\end{align}

\item[(iii)] Action-Type Conversion (AT-Conversion):  Player $i$ is able to choose $t_i$.
Define
\begin{align}\label{4.8}
\begin{array}{ccl}
c_i^{AC}(a,t_i):=\left[E(c_i(a)|t_i=t_i^1), E(c_i(a)|t_i=t_i^2),\cdots, E(c_i(a)|t_i=t_i^{\tau_i})\right],
\quad i=1,2,\cdots,n.
\end{array}
\end{align}
\end{itemize}

The corresponding complete information games are called Harsanyi Bayesian  game, Selten Bayesian  game, and Action-Type  Bayesian game of Bayesian game  $G$, respectively.
\end{dfn}

\begin{thm}\label{t4.7} The structure vectors of three conversions are as follows:
\begin{itemize}
\item[(i)]   H-Conversion: Denote by $$p=\left[Pr(t_1^1,\cdots,t_n^1),Pr(t_1^1,\cdots,t_n^2),\cdots,Pr(t_1^{\tau_1},\cdots,t_n^{\tau_n})\right]^T,$$
    then
\begin{align}\label{4.9}
V_i^H=V^{\bar{c}}_ip,\quad i=1,2,\cdots,n,
\end{align}
where $V_i^H$ is the  structure vector of $c_i^H(a).$
\item[(ii)]  S-Conversion: For a given $\bar{t}=(\bar{t}_1,\cdots,\bar{t}_n)$
\begin{align}\label{4.10}
V_i^S=V_i^{\bar{c}}W_{[\tau_i,\prod_{k=1}^{i-1}\tau_k]}\d_{\tau_i}^{\bar{t}_i}p_i(\bar{t}_i),\quad i=1,2,\cdots,n.
\end{align}
where $V_i^S$ is the  structure vector of $c_i^S(a).$
\item[(iii)]  AT-Conversion:
\begin{align}\label{4.11}
\begin{array}{ccl}
V_i^{AT}&=&V_i^{\bar{c}}W_{[\tau_i,\prod_{k=1}^{i-1}\tau_k]}\left[I_{\tau_i}*p_i\right],\qquad i=1,2,\cdots,n.
\end{array}
\end{align}
where $*$ is the Khatra-Rao product of matrices, and $V_i^{AT}$ is the  structure vector of $c_i^{AT}(a,t_i).$
\end{itemize}
\end{thm}

\noindent{\it Proof:}
\begin{itemize}
\item[(i)] According to Harsanyi Conversion,
\begin{align}\label{2.4}
\begin{array}{lllccc}
c_i^H(a)&=&\sum_{t\in T}Pr(t){\bar{c}}_i(a,t)\\
&~&\\
&=&\sum_{t\in T}Pr(t)V^{\bar{c}}_ita\\
&~&\\
&=&V^{\bar{c}}_ipa:=V_i^Ha.
\end{array}
\end{align}
\item[(ii)] According to Selten Conversion, player $i$ knows his type $t_i=\bar{t}_i$
\begin{align}\label{2.4}
\begin{array}{lllccc}
c_i^S(a)&=&E(c_i(a)|t_i=\bar{t}_i)\\
&~&\\
&=&\sum_{t_{-i}\in T_{-i}}p(t_{-i}|\bar{t}_i){\bar{c}}_i(a,\bar{t}_i,t_{-i})\\
&~&\\
&=&\sum_{t_{-i}\in T_{-i}}p(t_{-i}|\bar{t}_i)V^{\bar{c}}_iW_{[\tau_i,\prod_{k=1}^{i-1}\tau_k]}\bar{t}_it_{-i}a\\
&~&\\
&=&V_i^{\bar{c}}W_{[\tau_i,\prod_{k=1}^{i-1}\tau_k]}\d_{\tau_i}^{\bar{t}_i}p_i(\bar{t}_i)a\\
&~&\\
&=&V_i^Sa.
\end{array}
\end{align}
\item[(iii)] According to Action-Type Conversion,
\begin{align}\label{2.4}
\begin{array}{lllccc}
c_i^{AT}(a,t_i)&=&\left[E(c_i(a)|t_i=t_i^1), E(c_i(a)|t_i=t_i^2),\cdots, E(c_i(a)|t_i=t_i^{\tau_i})\right]\\
&~&\\
&=&\left[V_i^{\bar{c}}W_{[\tau_i,\prod_{k=1}^{i-1}\tau_k]}\d_{\tau_i}^{t_i^1}p_i(t_i^1)a,~\cdots,~ V_i^{\bar{c}}W_{[\tau_i,\prod_{k=1}^{i-1}\tau_k]}\d_{\tau_i}^{t_i^{\tau_i}}p_i(t_i^{\tau_i})a\right]\\
&~&\\
&=&V_i^{\bar{c}}W_{[\tau_i,\prod_{k=1}^{i-1}\tau_k]}\left[\d_{\tau_i}^{t_i^1}p_i(t_i^1),~\cdots,~ \d_{\tau_i}^{t_i^{\tau_i}}p_i(t_i^{\tau_i})\right]t_ia\\
&~&\\
&=&V_i^{\bar{c}}W_{[\tau_i,\prod_{k=1}^{i-1}\tau_k]}\left[I_{\tau_i}*p_i\right]t_ia\\
&~&\\
&=&V_i^{AT}t_ia.
\end{array}
\end{align}
\end{itemize}

\hfill $\Box$

\begin{dfn}\label{p4.8} Consider a  (finite) static Bayesian game $G=(N,T,A,c, Pr)$.

\begin{itemize}
\item[(i)] Harsanyi Bayesian game:

$a^*\in A$ is called a  Bayesian-Nash equilibrium for Harsanyi BG (H-BN-E), if $a^*$ satisfies
\begin{align}\label{4.12}
c_i^H(a_i^*,a_{-i}^*)\geq c_i^H(a_i,a_{-i}^*),\quad \forall a_i\in A_i,\; i=1,\cdots,n.
\end{align}

\item[(ii)] Selten Bayesian game:

$a^*\in A$ is called a  Bayesian-Nash equilibrium for Selten BG (S-BN-E) with respect to pre-assigned $\bar{t}$, if $a^*$ satisfies
\begin{align}\label{4.13}
c_i^S(a_i^*,a_{-i}^*)\geq c_i^S(a_i,a_{-i}^*),\quad \forall a_i\in A_i,\; i=1,\cdots,n.
\end{align}

\item[(iii)] Action-Type Bayesian game:

$(t^*, a^*)\in T\times A$ is called a  Bayesian-Nash equilibrium for Action-Type BG (AT-BN-E), if $(t^*, a^*)$ satisfies
\begin{align}\label{4.14}
\begin{array}{l}
c_i^{AT}(a_i^*,t_i^*,a_{-i}^*)\geq c_i^{AT}(a_i,t_i,a_{-i}^*),\quad \forall t_i\in T_i, \forall a_i\in A_i, \quad i=1,\cdots,n.
\end{array}
\end{align}
\end{itemize}

\end{dfn}

%

\begin{thm} \label{t4.10} Let O-BN-E be the original B-N equilibrium from Definition \ref{d4.3}. Compered with original BN, we have the following results:
\begin{itemize}
\item[(i)]
$$
a^*~\mbox{is an O-BN-E}\begin{array}{c}\Rightarrow\\\not\Leftarrow\end{array}a^*~\mbox{is an H-BN-E}.
$$
\item[(ii)]
$$
a^*~\mbox{is an O-BN-E}\begin{array}{c}\Rightarrow\\\Leftarrow\end{array}a^*~\mbox{is an S-BN-E for every}~ \bar{t}\in T.
$$
\item[(iii)]
$$
a^*~\mbox{is an O-BN-E}\begin{array}{c}\not\Rightarrow\\\Leftarrow\end{array}(t^*,a^*)~\mbox{is an AT-BN-E}.
$$
\end{itemize}
\end{thm}

\begin{exa}\label{e4.11} Recall the Bayesian game in Example \ref{e3.3}.
\begin{itemize}
  \item[(i)] Harsanyi game:

  It is easy to calculate that
$$
\begin{array}{l}
V^H_1=[-\infty,-\infty,-\infty,1.6,-\infty,-0.6,-\infty,-\infty,-\infty],\\
V^H_2=[-\infty,-\infty,-\infty,-0.1,-\infty,-2.1,-\infty,-\infty,-\infty],\\
\end{array}
$$
Put them into bi-matrix form:

\begin{table}[!htb]
\centering
\caption{Harsanyi Expected Payoff Bi-Matrix}
\vskip .25\baselineskip
\begin{tabular}{|c||c|c|c|}
\hline
$a_1\backslash a_2$&$\d_3^1$&$\d_3^2$&$\d_3^3$\\
\hline
\hline
$\d_3^1$&$-\infty,~-\infty$&$-\infty,~-\infty$&$-\infty,~-\infty$\\
\hline
$\d_3^2$&{\color{red}$1.6,~-0.1$}&$-\infty,~-\infty$&$-0.6,~-2.1$\\
\hline
$\d_3^3$&$-\infty,~-\infty$&$-\infty,~-\infty$&$-\infty,~-\infty$\\
\hline
\end{tabular}
\end{table}

Hence $(a_1^2,a_2^1)=(\d_3^2,\d_3^1)$ is an H-BN-E.
  \item[(ii)] Action-Type game:

  We can also calculate that
$$
\begin{array}{l}
\begin{array}{ccccccccccccc}
V^{AT}_1&=&[1.3&-\infty&1&1&-\infty&-0.5&-\infty&-\infty&-\infty\\
~&~~&-\infty&-\infty&-\infty&2.3&-\infty&0&2&-\infty&0.3],\\
\end{array}\\
~\\
\begin{array}{ccccccccccccc}
V^{AT}_2&=&[-\infty&-\infty&-\infty&0.8&2.2&-1.8&-\infty&-\infty&-\infty\\
~&~~&-\infty&-\infty&-\infty&0.4&-\infty&-2.4&-\infty&-\infty&-\infty].
\end{array}
\end{array}
$$

Put them into bi-matrix form (Table \ref{tb.10} and Table \ref{tb.11}):
\begin{table}[!htb]
\centering
\caption{TH: Expected Payoff Bi-Matrix}\label{tb.10}
\vskip .25\baselineskip
\begin{tabular}{|c||c|c|c|}
\hline
$t_1a_1\backslash t_2a_2$&$t_2^1a_2^1$&$t_2^2a_2^1$&$t_2^1a_2^2$\\
\hline
\hline
$t_1^1a_1^1$&$1.3,~-\infty$&$1.3.~-\infty$&$-2.5,~-\infty$\\
\hline
$t_1^2a_1^1$&$-\infty,~-\infty$&$-\infty,~-\infty$&$-\infty,~-\infty$\\
\hline
$t_1^1a_1^2$&$1,~0.8$&$1,~0.4$&$-\infty,~2.2$\\
\hline
$t_1^2a_1^2$&$2.3,~0.8$&$2.3,~0.4$&$-\infty,~2.2$\\
\hline
$t_1^1a_1^3$&$-\infty,~-\infty$&$-\infty,~-\infty$&$-\infty,~-\infty$\\
\hline
$t_1^2a_1^3$&$2,~-\infty$&$2,~-\infty$&$-\infty,~-\infty$\\
\hline
\end{tabular}
\end{table}

\begin{table}[!htb]
\centering
\caption{A-T- Expected Payoff Bi-Matrix(cont'd)}\label{tb.11}
\vskip .25\baselineskip
\begin{tabular}{|c||c|c|c|}
\hline
$t_1a_1\backslash t_2a_2$&$t_2^2a_2^2$&$t_2^1a_2^3$&$t_2^2a_2^3$\\
\hline
\hline
$t_1^1a_1^1$&$-2.5,~-\infty$&$1,~-\infty$&$1,~-\infty$\\
\hline
$t_1^2a_1^1$&$-\infty,~-\infty$&$-\infty,~-\infty$&$-\infty,~-\infty$\\
\hline
$t_1^1a_1^2$&$-\infty,~-\infty$&$-0.5,~-1.8$&$-0.5,~-2.4$\\
\hline
$t_1^2a_1^2$&$-\infty,~-\infty$&$0,~-1.8$&$0,~-2.4$\\
\hline
$t_1^1a_1^3$&$-\infty,~-\infty$&$-\infty,~-\infty$&$-\infty,~-\infty$\\
\hline
$t_1^2a_1^3$&$-\infty,~-\infty$&$0.3,~-\infty$&$0.3,~-\infty$\\
\hline
\end{tabular}
\end{table}

It is easy to verify that there is no AT-BN-E.

\end{itemize}

\end{exa}

\section{Bayesian Potential Game}

\begin{dfn}\label{d5.1} \cite{heu96} Consider a  (finite) static Bayesian game $G=(N,T,A,c, Pr)$.
It is called a TN  weighted potential game, if there exists a function $F:T\times A\ra \R$ such that for any $t\in T$
\begin{align}\label{5.1}
\begin{array}{l}
c_i(a_i',a_{-i},t)-c_i(a_i,a_{-i},t)=w_i\left(F(a_i',a_{-i},t)-F(a_i,a_{-i},t)\right),\\
~~~~~a'_i,a_i\in A_i, \;a_{-i}\in A_{-i},\;t\in T,
\end{array}
\end{align}
where $w_i>0$. When $w_i=1$, it is a  TN potential game, $\forall i$.
\end{dfn}

\begin{dfn}\label{d5.2} Consider a  (finite) static Bayesian game $G=(N,T,A,c, Pr)$. It is called a TH weighted potential game, if there exists a function $F:T\times A\ra \R$ such that for any $t_i,t_i'\in T_i,~a_i,a_i'\in A_i$
\begin{align}\label{5.2}
\begin{array}{l}
c_i(a_i',a_{-i},t_i',t_{-i})-c_i(a_i,a_{-i},t_i,t_{-i})=w_i\left(F(a_i',a_{-i},t_i',t_{-i})-F(a_i,a_{-i},t_i,t_{-i})\right),\\
~~~~~ ~~~~~~~~~~~~~~~~~~~~~~\forall \;a_{-i}\in A_{-i}, t_{-i}\in T_{-i}.
\end{array}
\end{align}
where $w_i>0$. When $w_i=1$, it is  a TH potential game, $\forall i$.
\end{dfn}

The following proposition comes from definition immediately.

\begin{prp}\label{p5.3}
$$
\mbox{TH potential game} ~\Rightarrow \mbox{TN potential game}.
$$
\end{prp}

\begin{rem}\label{r5.5} \cite{gffp} It is obvious that a Bayesian game $G=(N,T,A,c, Pr)$ is TN-potential, if and only if, under each type it is potential. Hence, we can use potential equation [6] to check whether a game is potential and to construct potential function. Please refer to Theorem \ref{t2.3.5}.
\end{rem}

\begin{exa}\label{e5.4} Consider a finite BG
\begin{align}\label{5.3}
G=(N,T,A,c,p),
\end{align}
where
$$
\begin{array}{ccl}
N&=&\{1,2\};\\
T&=&\{T_1,T_2\},\;\mbox{with}\\
T_1&=&\{t_1^1,t_1^2,t_1^3\},\\
T_2&=&\{t_2^1,t_2^2\};\\
A&=&\{A_1,A_2\},\;\mbox{with}\\
A_1&=&\{a_1^1,a_1^2\},\\
A_2&=&\{a_2^1,a_2^2,a_2^3\},\\
\end{array}
$$

$T$ is as follows:
\begin{table}[!htb]
\centering
\caption{Distribution}\label{tb.5.1}
\vskip .25\baselineskip
\begin{tabular}{|c||c|c|}
\hline
$t_1\backslash t_2$&$t_2^1$&$t_2^2$\\
\hline
\hline
$t_1^1$&$0.1$&$0.15$\\
\hline
$t_1^2$&$0.15$&$0.2$\\
\hline
$t_1^3$&$0.3$&$0.1$\\
\hline
\end{tabular}
\end{table}

\vskip 2mm

Then
$$
p_1({t_1^1})=\begin{bmatrix}
Pr(t_2^1|t_1^1)\\
Pr(t_2^2|t_1^1)\\
\end{bmatrix}=\begin{bmatrix}
0.4\\
0.6\\
\end{bmatrix}.
$$

Similarly, we have

$$
\begin{array}{ll}
p_1({t_1^2})=(3/7,4/7)^T,&p_1({t_1^3})=(0.75,0.25)^T,\\
p_2({t_2^1})=(2/11,3/11,6/11)^T,&p_2({t_2^2})=(1/3,4/9,2/9)^T.\\
\end{array}
$$

The payoff functions are as follows.
\begin{table}[!htb]
\centering
\caption{Payment Function for $t_1^1-t_2^1$}\label{tb.5.2}
\vskip .25\baselineskip
\begin{tabular}{|c||c|c|c|}
\hline
$c_1\backslash c_2$&$a_2^1$&$a_2^2$&$a_2^3$\\
\hline
\hline
$a_1^1$&$5,~0$&$2,~2$&$0,~1$\\
\hline
$a_1^2$&$2,~-1$&$-1,~1$&$1,~4$\\
\hline
\end{tabular}
\end{table}

\vskip -2mm

\begin{table}[!htb]
\centering
\caption{Payment Function for $t_1^1-t_2^2$}\label{tb.5.3}
\vskip .25\baselineskip
\begin{tabular}{|c||c|c|c|}
\hline
$c_1\backslash c_2$&$a_2^1$&$a_2^2$&$a_2^3$\\
\hline
\hline
$a_1^1$&$3,~0$&$2,~3$&$1,~1$\\
\hline
$a_1^2$&$1,~-2$&$0,~1$&$2,~2$\\
\hline
\end{tabular}
\end{table}

\vskip -2mm

\begin{table}[!htb]
\centering
\caption{Payment Function for $t_1^2-t_2^1$}\label{tb.5.4}
\vskip .25\baselineskip
\begin{tabular}{|c||c|c|c|}
\hline
$c_1\backslash c_2$&$a_2^1$&$a_2^2$&$a_2^3$\\
\hline
\hline
$a_1^1$&$2,~4$&$0,~-1$&$5,~5$\\
\hline
$a_1^2$&$1,~1$&$3,~0$&$1,~-1$\\
\hline
\end{tabular}
\end{table}

\vskip -2mm

\begin{table}[!htb]
\centering
\caption{Payment Function for $t_1^2-t_2^2$}\label{tb.5.5}
\vskip .25\baselineskip
\begin{tabular}{|c||c|c|c|}
\hline
$c_1\backslash c_2$&$a_2^1$&$a_2^2$&$a_2^3$\\
\hline
\hline
$a_1^1$&$1,~0$&$-1,~2$&$1,~3$\\
\hline
$a_1^2$&$3,~-1$&$0,~0$&$2,~1$\\
\hline
\end{tabular}
\end{table}

\vskip -2mm

\begin{table}[!htb]
\centering
\caption{Payment Function for $t_1^3-t_2^1$}\label{tb.5.6}
\vskip .25\baselineskip
\begin{tabular}{|c||c|c|c|}
\hline
$c_1\backslash c_2$&$a_2^1$&$a_2^2$&$a_2^3$\\
\hline
\hline
$a_1^1$&$1,~1$&$2,~2$&$0,~3$\\
\hline
$a_1^2$&$4,~2$&$5,~3$&$2,~3$\\
\hline
\end{tabular}
\end{table}

\vskip -2mm

\begin{table}[!htb]
\centering
\caption{Payment Function for $t_1^3-t_2^2$}\label{tb.5.7}
\vskip .25\baselineskip
\begin{tabular}{|c||c|c|c|}
\hline
$c_1\backslash c_2$&$a_2^1$&$a_2^2$&$a_2^3$\\
\hline
\hline
$a_1^1$&$-2,~1$&$1,~2$&$0,~0$\\
\hline
$a_1^2$&$-4,~1$&$-5,~-2$&$-2,~0$\\
\hline
\end{tabular}
\end{table}

The potential function is obtained as
\begin{table}[H]
\centering
\caption{Potential Function for $t_1^1-t_2^1$}\label{tb.5.8}
\vskip .25\baselineskip
\begin{tabular}{|c||c|c|c|}
\hline
$F$&$a_2^1$&$a_2^2$&$a_2^3$\\
\hline
\hline
$a_1^1$&$1$&$3$&$2$\\
\hline
$a_1^2$&$-2$&$0$&$3$\\
\hline
\end{tabular}
\end{table}

\vskip 2mm

\begin{table}[H]
\centering
\caption{Potential Function for $t_1^1-t_2^2$}\label{tb.5.9}
\vskip .25\baselineskip
\begin{tabular}{|c||c|c|c|}
\hline
$F$&$a_2^1$&$a_2^2$&$a_2^3$\\
\hline
\hline
$a_1^1$&$1$&$4$&$2$\\
\hline
$a_1^2$&$-1$&$2$&$3$\\
\hline
\end{tabular}
\end{table}

\vskip 2mm

\begin{table}[!htb]
\centering
\caption{Potential Function for $t_1^2-t_2^1$}\label{tb.5.10}
\vskip .25\baselineskip
\begin{tabular}{|c||c|c|c|}
\hline
$F$&$a_2^1$&$a_2^2$&$a_2^3$\\
\hline
\hline
$a_1^1$&$3$&$-2$&$4$\\
\hline
$a_1^2$&$2$&$1$&$0$\\
\hline
\end{tabular}
\end{table}

\vskip 2mm

\begin{table}[!htb]
\centering
\caption{Potential Function for $t_1^2-t_2^2$}\label{tb.5.11}
\vskip .25\baselineskip
\begin{tabular}{|c||c|c|c|}
\hline
$F$&$a_2^1$&$a_2^2$&$a_2^3$\\
\hline
\hline
$a_1^1$&$-1$&$1$&$2$\\
\hline
$a_1^2$&$1$&$2$&$3$\\
\hline
\end{tabular}
\end{table}

\vskip 2mm

\begin{table}[!htb]
\centering
\caption{Potential Function for $t_1^3-t_2^1$}\label{tb.5.12}
\vskip .25\baselineskip
\begin{tabular}{|c||c|c|c|}
\hline
$F$&$a_2^1$&$a_2^2$&$a_2^3$\\
\hline
\hline
$a_1^1$&$-1$&$0$&$1$\\
\hline
$a_1^2$&$2$&$3$&$3$\\
\hline
\end{tabular}
\end{table}

\vskip 2mm

\begin{table}[!htb]
\centering
\caption{Potential Function for $t_1^3-t_2^2$}\label{tb.5.13}
\vskip .25\baselineskip
\begin{tabular}{|c||c|c|c|}
\hline
$F$&$a_2^1$&$a_2^2$&$a_2^3$\\
\hline
\hline
$a_1^1$&$4$&$5$&$3$\\
\hline
$a_1^2$&$2$&$-1$&$1$\\
\hline
\end{tabular}
\end{table}

Finally, we have
$$
\begin{array}{ccl}
V^c_1&=&[5,2,0,2,-1,1,3,2,1,1,0,2,2,0,5,1,3,1,\\
~&~&1,-1,1,3,0,2,1,2,0,4,5,2,-2,1,0,-4,-5,-2],\\
V^c_2&=&[0,2,1,-1,1,4,0,3,1,-2,1,2,4,-1,5,1,0,\\
~&~&-1,0,2,3,-1,0,1,1,2,3,2,3,3,1,2,0,1,-2,0],\\
V^c_F&=&[2,4,3,-1,1,4,0,3,1,-1,1,2,2,-3,3,1,0,-1,\\
~&~&-3,-1,0,-1,0,1,-1,0,1,2,3,3,3,4,2,1,-2,0].\\
\end{array}
$$

It is easy to verify that $F$ is a TN potential function.
\end{exa}

\begin{dfn}\label{p4.8} Consider a  finite static Bayesian game $G=(N,T,A,c, Pr)$.
\begin{itemize}
\item[(i)] It is called a Harsanyi weighted potential game, if there exists a function $Q^H:A\ra \R$ such that for any $a_i,a_i'\in A_i$
\begin{align}\label{5.2}
\begin{array}{l}
c_i^H(a_i',a_{-i})-c_i^H(a_i,a_{-i})=w_i\left(Q^H(a_i',a_{-i})-Q^H(a_i,a_{-i})\right),~~~\forall \;a_{-i}\in A_{-i}.
\end{array}
\end{align}
where $w_i>0$. When $w_i=1$, $\forall i$, it is  a Harsanyi potential game.

\item[(ii)] It is called a Selten weighted potential game for the pre-assigned $t\in T$, if there exists a function $Q^{S}:A\ra \R$ such that for any $a_i,a_i'\in A_i$
\begin{align}\label{5.2}
\begin{array}{l}
c_i^S(a_i',a_{-i})-c_i^S(a_i,a_{-i})=w_i\left(Q^S(a_i',a_{-i})-Q^S(a_i,a_{-i})\right),~~~\forall \;a_{-i}\in A_{-i}.
\end{array}
\end{align}
where $w_i>0$. When $w_i=1$, $\forall i$, it is  a Selten potential game.

\item[(iii)]
It is called a Action-Type weighted potential game, if there exists a function $Q^{AT}:T\times A\ra \R$ such that for any $a_i,a_i'\in A_i$
\begin{align}\label{eq68}
\begin{array}{l}
c_i^{AT}(a_i',t_i',a_{-i})-c_i^{AT}(a_i,t_i,a_{-i})=w_i\left(Q^{AT}(a_i',t_i',a_{-i},t_{-i})-Q^{AT}(a_i,t_i,a_{-i},t_{-i})\right),~\forall \;a_{-i}\in A_{-i}.
\end{array}
\end{align}
where $w_i>0$. When $w_i=1$, $\forall i$, it is  a Action-Type potential game.

\end{itemize}

\end{dfn}

The following proposition is obvious.
\begin{prp}
Consider a  finite static Bayesian game $G=(N,T,A,c, Pr)$.
\begin{itemize}
\item[(i)] If $G$ is a Harsanyi  potential game, then $G$ has at least one H-BN-E (the potential maximizer).
\item[(ii)] If $G$ is a Selten  potential game, then $G$ has at least one S-BN-E (the potential maximizer).
\item[(iii)] If $G$ is a Action-Type  potential game, then $G$ has at least one AT-BN-E (the potential maximizer).
\end{itemize}
\end{prp}

\begin{prp}\label{r5.5}  Consider a  finite  Bayesian game $G=(N,T,A,c, Pr)$. If $G$ is a TN or TH potential game with potential function $F(a,t)$, then the corresponding  Harsanyi Bayesian game  is a Harsanyi  potential game with potential function
\begin{align}\label{5.1}
\begin{array}{lllcc}
Q^H(a)&=&\sum_{t\in T}Pr(t)F(a,t)\\
&=&V_F^cpta\\
&:=&V^Ha,\\
\end{array}
\end{align}
where $V_F^c$ is the structure vector of $F(a,t).$
\end{prp}

\begin{exa}\label{e5.7}  Consider Example \ref{e5.4} again.
\begin{itemize}
\item[(i)] Harsanli Bayesian game:

The payoff vectors of corresponding Harsanli BG are
$$
\begin{array}{ccl}
V^H_1&=&V^c_1p=[1.55,1,1.1,1.9,1.35,1.35],\\
V^H_2&=&V^c_2p=[1,1.7,2.5,0.25,0.95,1.65].\\
V^c_F&=&[2,4,3,-1,1,4,0,3,1,-1,1,2,2,-3,3,1,0,-1,\\
~&~&-3,-1,0,-1,0,1,-1,0,1,2,3,3,3,4,2,1,-2,0].\\
\end{array}
$$
Back to matrix form:
\begin{table}[!htb]
\centering
\caption{Harsanli-Bayesian }\label{tb.5.14}
\vskip .25\baselineskip
\begin{tabular}{|c||c|c|c|}
\hline
$c_1^H\backslash c_2^H$&$a_2^1$&$a_2^2$&$a_2^3$\\
\hline
\hline
$a_1^1$&$1.55,~1$&$1,~1.7$&$1.1,~2.5$\\
\hline
$a_1^2$&$1.9,~0.25$&$1.35,~0.95$&$\underline{1.35},~\underline{1.65}$\\
\hline
\end{tabular}
\end{table}

\vskip 5mm

It is easy to verify that the Harsanli Bayesian game is  potential with $V^H$ as its potential function, where
$$V^H=V^c_Fp=[-0.1,0.6,1.4,0.25,0.95,1.65].$$
Back to matrix form:
\begin{table}[!htb]
\centering
\caption{Harsanli Potential Function $V^H$}\label{tb.5.15}
\vskip .25\baselineskip
\begin{tabular}{|c||c|c|c|}
\hline
$Q_H$&$a_2^1$&$a_2^2$&$a_2^3$\\
\hline
\hline
$a_1^1$&$-0.1$&$0.6$&$1.4$\\
\hline
$a_1^2$&$0.25$&$0.95$&$1.65$\\
\hline
\end{tabular}
\end{table}
\item[(ii)] Selten Bayesian game:

$$
V^S_{i}=V^c_iW_{[\tau_i,\prod_{j=1}^{i-1}\tau_j]}\d_{\tau_i}^\theta p_i({t_i^\theta}),\quad 1\leq \theta\leq \tau_i,\;i=1,\cdots,n.
$$

Assume $t_1^{\theta}=t_1^1$, $t_2^{\theta}=t_2^1$. Then we have

$$
\begin{array}{ccl}
V^S_{1}&=&[3.8,2,0.6,1.4,-0.4,1.6],\\
V^S_{2}&=&[2.1818,    0.6364,    3.1818,   -0.1818,    0.4545,    0.7273].
\end{array}
$$

Back to matrix form

\vskip 2mm

\begin{table}[!htb]
\centering
\caption{Selten-Bayesian}\label{tb.5.16}
\vskip .25\baselineskip
\begin{tabular}{|c||c|c|c|}
\hline
$P_1\backslash P_2$&$a_2^1$&$a_2^2$&$a_2^3$\\
\hline
\hline
$a_1^1$&$3.8,~2.1818$&$2,~0.6364$&$0.6,~3.1818$\\
\hline
$a_1^2$&$1.4, -0.1818$&$-0.4,~0.4545$&$1.6,~0.7273$\\
\hline
\end{tabular}
\end{table}

\vskip 2mm

Assume $t_1^{\theta}=t_1^1$, $t_2^{\theta}=t_2^2$. Then we have

$$
\begin{array}{ccl}
V^S_{11}&=&[3.8,2,0.6,1.4,-0.4,1.6],\\
V^S_{22}&=&[ 0.2222,    2.3333,    1.6667,   -0.8889,   -0.1111,    1.1111].
\end{array}
$$

Back to matrix form

\vskip 2mm

\begin{table}[!htb]
\centering
\caption{Selten-Bayesian}\label{tb.5.17}
\vskip .25\baselineskip
\begin{tabular}{|c||c|c|c|}
\hline
$P_1\backslash P_2$&$a_2^1$&$a_2^2$&$a_2^3$\\
\hline
\hline
$a_1^1$&$3.8,~0.2222$&$2,~2.3333$&$0.6,~1.6667$\\
\hline
$a_1^2$&$1.4, ~-0.8889$&$-0.4,~-0.1111$&$1.6,~1.1111$\\
\hline
\end{tabular}
\end{table}

Assume $t_1^{\theta}=t_1^2$, $t_2^{\theta}=t_2^1$. Then we have

$$
\begin{array}{ccl}
V^S_{12}&=&[1.8571,    2,    0.4286,    2.7143,    2.8571,    2],\\
V^S_{21}&=&[2.1818,    0.6364,    3.1818,   -0.1818,    0.4545,    0.7273].
\end{array}
$$

Back to matrix form

\vskip 2mm

\begin{table}[!htb]
\centering
\caption{Selten-Bayesian}\label{tb.5.18}
\vskip .25\baselineskip
\begin{tabular}{|c||c|c|c|}
\hline
$P_1\backslash P_2$&$a_2^1$&$a_2^2$&$a_2^3$\\
\hline
\hline
$a_1^1$&$1.8571,~2.1818$&$2,~0.6364$&$0.4286,~3.1818$\\
\hline
$a_1^2$&$2.7143, -0.1818$&$2.8571,~0.4545$&$2,~0.7273$\\
\hline
\end{tabular}
\end{table}

\vskip 2mm

\begin{table}[!htb]
\centering
\caption{Selten-Bayesian}\label{tb.5.19}
\vskip .25\baselineskip
\begin{tabular}{|c||c|c|c|}
\hline
$P_1\backslash P_2$&$a_2^1$&$a_2^2$&$a_2^3$\\
\hline
\hline
$a_1^1$&$1.8571,~2$&$2,~2.3333$&$0.4286,~1.6667$\\
\hline
$a_1^2$&$2.7143, ~-0.8889$&$2.8571,~-0.1111$&$2,~1.1111$\\
\hline
\end{tabular}
\end{table}

\vskip 2mm

\begin{table}[!htb]
\centering
\caption{Selten-Bayesian}\label{tb.5.20}
\vskip .25\baselineskip
\begin{tabular}{|c||c|c|c|}
\hline
$P_1\backslash P_2$&$a_2^1$&$a_2^2$&$a_2^3$\\
\hline
\hline
$a_1^1$&$1,~2.1818$&$0.25,~0.6364$&$3.75,~3.1818$\\
\hline
$a_1^2$&$-0.25, -0.1818$&$1,~0.4545$&$0.25,~0.7273$\\
\hline
\end{tabular}
\end{table}

Assume $t_1^{\theta}=t_1^3$, $t_2^{\theta}=t_2^2$. Then we have

$$
\begin{array}{ccl}
V^S_{13}&=&[1,0.25,    3.75,   -0.25,    1,    0.25],\\
V^S_{22}&=&[ 0.2222,    2.3333,    1.6667,   -0.8889,   -0.1111,    1.1111].
\end{array}
$$

Back to matrix form

\vskip 2mm

\begin{table}[!htb]
\centering
\caption{Selten-Bayesian}\label{tb.5.21}
\vskip .25\baselineskip
\begin{tabular}{|c||c|c|c|}
\hline
$P_1\backslash P_2$&$a_2^1$&$a_2^2$&$a_2^3$\\
\hline
\hline
$a_1^1$&$1,~2$&$0.25,~2.3333$&$3.75,~1.6667$\\
\hline
$a_1^2$&$-0.25, ~-0.8889$&$1,~-0.1111$&$0.25,~1.1111$\\
\hline
\end{tabular}
\end{table}

It is easy to verify that the Selten game is not potential for each type.
\item[(iii)] Action-Type game:

$$
V^{AT}_i=\left[V^S_{i,1},V^S_{i,2},\cdots,V^S_{i,s_i}\right],\quad i=1,2,\cdots,n.
$$

$$
\begin{array}{lllllll}
V^{AT}_1=&~&~&~&~&~&~\\
~&[3.8&2&0.6&1.4&-0.4&1.6\\
~&1.8571&    2&    0.4286& 2.7143&    2.8571&    2\\
~&  1&    0.25&    3.75&   -0.25&    1&    0.25],\\
V^{AT}_2=&~&~&~&~&~&~\\
~&[2.1818&    0.6364&    3.1818&   -0.1818&    0.4545&    0.7273  \\
~&0.2222&    2.3333&    1.6667&   -0.8889&   -0.1111&    1.1111].
\end{array}
$$

\vskip 2mm

Back to matrix form
\begin{table}[!htb]
\centering
\caption{Action-Type-Bayesian}\label{tb.5.22}
\vskip .25\baselineskip
\begin{tabular}{|c||c|c|c|}
\hline
$P_1\backslash P_2$&$t_2^1a_2^1$&$t_2^2a_2^1$&$t_2^1a_2^2$\\
\hline
\hline
$t_1^1a_1^1$&$3.8,~    2.1818$&$3.8,~    0.6364$&$    2,~    3.1818$\\
\hline
$t_1^2a_1^1$&$1.4,~    2.1818$&$1.4,~0.6364$&$   -0.4,~3.1818$\\
\hline
$t_1^3a_1^1$&$ 1.8571,    2.1818$&$    1.8571,~    0.6364$&$    2,~3.1818 $\\
\hline
$t_1^1a_1^2$&$ 2.7143,    0.2222$&$    2.7143,~    2.3333$&$    2.8571,~    1.6667 $\\
\hline
$t_1^2a_1^2$&$1,    0.2222$&$    1,~    2.3333$&$    0.25,~    1.6667   $\\
\hline
$t_1^3a_1^3$&$ -0.25,~    0.2222$&$   -0.25,~    2.3333$&$    1,    1.6667  $\\
\hline
\end{tabular}
\end{table}

\vskip 2mm

\begin{table}[!htb]
\centering
\caption{Action-Type-Bayesian}\label{tb.5.23}
\vskip .25\baselineskip
\begin{tabular}{|c||c|c|c|}
\hline
$P_1\backslash P_2$&$t_2^2a_2^2$&$t_2^1a_2^3$&$t_2^2a_2^3$\\
\hline
\hline
$t_1^1a_1^1$&$ 2,   -0.1818$&$0.6,    0.4545$&$    0.6,    0.7273$\\
\hline
$t_1^2a_1^1$&$  -0.4,~   -0.1818$&$    1.6,    0.4545$&$    1.6,    0.7273$\\
\hline
$t_1^3a_1^1$&$  2,~  -0.1818$&$    0.4286,    0.4545$&$    0.4286,~    0.7273$\\
\hline
$t_1^1a_1^2$&$ 2.8571,~   -0.8889$&$    2,~   -0.1111$&$    2,~    1.1111 $\\
\hline
$t_1^2a_1^2$&$       0.25,   -0.8889$&$    3.75,~   -0.1111$&$    3.75,~    1.1111$\\
\hline
$t_1^3a_1^3$&$     1,~   -0.8889$&$    0.25,~   -0.1111$&$    0.25,~    1.1111 $\\
\hline
\end{tabular}
\end{table}

\vskip 2mm

The Action-Type game is not potential too.

\end{itemize}
\end{exa}

According to the definition,  a Bayesian game  is a Harsanli potential game, if and only if, the original game is potential for each type. Hence, we can use potential equation (\ref{2.3.8}) to check whether a Bayesian game is Harsanli potential and to construct potential function. But the check for Selten potential game and Action-Type potential game are not obvious.

In the following, we provide a method to verify whether a BG is Selten potential or  Action-Type potential.
Denote by
$$
\phi_i=\otimes_{j=1}^n\gamma_j,
$$
where
$$
\gamma_j=
\begin{cases}
I_{\tau_i},\quad j= i,\\
{\bf 1}_{\tau_j}^T,\quad j\neq i.
\end{cases}
$$

And let $W_G^S$ and $W_G^{AT}$ be  equivalent vectors of Selten Bayesian game and Action-Type Bayesian game respectively
\begin{align}\label{3.13}
\begin{array}{lcc}
W_G^S&:=&[V_2^S-V_1^S,V_3^S-V_1^S,\cdots,V_n^S-V_1^S],~~~~~~~\\
~&~&~\\
W_G^{AT}&:=&[V_2^{AC}\phi_2-V_1^{AC}\phi_1,\cdots,V_n^{AC}\phi_n-V_1^{AC}\phi_1].
\end{array}
\end{align}
It is easy to verify that
$$W_G^S=V_G\Gamma^S_{\bar{t}},~~W_G^{AT}=V_G\Gamma^{AT},$$
where $V_G=[V_1^{\bar{c}},\cdots,V_n^{\bar{c}}],~\Gamma^S_{t_i}=W_{[\tau_i,\prod_{k=1}^{i-1}\tau_k]}\d_{\tau_i}^{\bar{t}_i}p_i(\bar{t}_i),~~~\Gamma^{AT}_i=W_{[\tau_i,\prod_{k=1}^{i-1}\tau_k]}\left[I_{\tau_i}*p_i\right]\phi_i,~\forall i,$
\begin{align}\label{eq71}
\Gamma^S_{\bar{t}}=\begin{bmatrix}
-\Gamma^S_{t_1}&-\Gamma^S_{t_1}&\cdots&-\Gamma^S_{t_1}\\
\Gamma^S_{t_2}&0&\cdots&0\\
0&\Gamma^S_{t_3}&\cdots&0\\
\vdots&~&~&~\\
0&0&\cdots&-\Gamma^S_{t_n}
\end{bmatrix},~~
\Gamma^{AT}=\begin{bmatrix}
-\Gamma^{AT}_1&-\Gamma^{AT}_1&\cdots&-\Gamma^{AT}_1\\
\Gamma^{AT}_2&0&\cdots&0\\
0&\Gamma^{AT}_3&\cdots&0\\
\vdots&~&~&~\\
0&0&\cdots&-\Gamma^{AT}_n
\end{bmatrix}.
\end{align}

\vskip 2mm
\begin{thm}\label{thm5.10}  Consider a  finite  Bayesian game $G=(N,T,A,c, Pr)$.
\begin{itemize}
\item[(i)]
$G$ is a Selten potential game, if and only if, the following linear equation has a solution
\begin{align}\label{eq72}
\Psi\xi=b^S,
\end{align}
where $\xi=(\xi_1^T,\cdots,\xi_n^T)^T,~b^S=(W_G^S)^T$, and $\Psi$ is defined by (\ref{2.3.9}).
Moreover, when a solution exists,   the structure vector of Selten potential function $Q^S(a)$ is
\begin{align}\label{2.3.10}
\begin{array}{llcc}
V^S_Q&=&V^S_1-\xi_1^TE_1~~~~~~~~~~~~\\
&=&V^c_1\d_{\tau_1}^{\bar{t}_1}p_1(\bar{t}_1)-\xi_1^TE_1.
\end{array}
\end{align}
\item[(ii)] $G$ is a Action-Type potential game, if and only if, the following linear equation has a solution
\begin{align}\label{eq74.1}
\Psi^{AT}\xi=b^{AT},
\end{align}
where $\xi=(\xi_1^T,\cdots,\xi_n^T)^T,~b^{AT}=(W_G^{AT})^T,~\varphi_i=I_{\theta_i}\otimes {\bf 1}^T_{\tau_i}\otimes I_{\vartheta_i}\otimes {\bf 1}^T_{r_i},  \vartheta_i=\prod_{j=i+1}^{n}\tau_i\prod_{l=1}^{i-1}r_l$,
\begin{align}
\theta_i=\prod_{j=1}^{i-1}\tau_i,~~
\Psi^{AT}=\begin{bmatrix}
-\varphi_1^T&\varphi_2^T&0&\cdots&0\\
-\varphi_1^T&0&\varphi_3^T&\cdots&0\\
\vdots&~&~&\ddots&~\\
-\varphi_1^T&0&0&\cdots&\varphi_n^T
\end{bmatrix}.
\end{align}
Moreover, when a solution exists,  the structure vector of Action-Type potential function $Q^{AT}(a,t)$ is
\begin{align}\label{2.3.10}
\begin{array}{llcc}
V^{AT}_Q&=&V^{AT}_1-\xi_1^T\varphi_1~~~~~~~~~~\\
&=&V^c_1\left[I_{\tau_i}*p_i\right]\phi_1-\xi_1^T\varphi_1.
\end{array}
\end{align}
\end{itemize}
\end{thm}

\noindent{\it Proof:}
The proof of (i) is similar with Theorem \ref{t2.3.5}, so we only prove (ii). According to (\ref{eq68}), there exist functions $d_i(a_{-i},t_{-i},\hat{a}_i,\hat{t}_i)$ such that
\begin{align}\label{eq76}
c_i^{AT}(a,t_i)=Q^{AT}(a,t)+d_i(a_{-i},t_{-i},\hat{a}_i,\hat{t}_i),~\forall i,
\end{align}
where the symbol ``hat" means that $d_i(a_{-i},t_{-i},\hat{a}_i,\hat{t}_i)$ is independent of $a_i$ and $t_i.$ Taking its vector form, (\ref{eq76}) becomes
\begin{align}\label{eq77}
\begin{array}{llcc}
~~~~V_i^{AT}t_ia&=&V_Q^{AT}ta+V^d_it_{-i}a_{-i}\\
\Leftrightarrow V_i^{AT}\phi_ita&=&V_Q^{AT}ta+V^d_i\varphi_ita.
\end{array}
\end{align}
Equation (\ref{eq77}) is equivalent to
\begin{align}\label{eq78}
V_Q^{AT}=V_i^{AT}\phi_i -V^d_i\varphi_i.
\end{align}
Then
\begin{align}\label{eq79}
V_i^{AT}\phi_i -V_1^{AT}\phi_1=V^d_i\varphi_i -V^d_1\varphi_1,~i=2,3,\cdots,n.
\end{align}
Equation (\ref{eq79}) is equivalent to the following linear equations
$$\Psi^{AT}\xi=b^{AT},$$
where $\xi_i=(V_i^d)^T,~i=1,2,\cdots,n.$
\hfill $\Box$

\begin{exa}\label{e5.11} Consider a finite Bayesian game $G=(N,T,A,c,Pr)$, where
$$
\begin{array}{ccl}
N&=&\{1,2\};\\
T&=&\{T_1,T_2\},\;\mbox{with}\\
T_1&=&\{t_1^1,t_1^2\},\\
T_2&=&\{t_2^1,t_2^2\};\\
A&=&\{A_1,A_2\},\;\mbox{with}\\
A_1&=&\{a_1^1,a_1^2\},\\
A_2&=&\{a_2^1,a_2^2\},\\
\end{array}
$$

The distribution $Pr$ is as follows:
\begin{table}[!htb]
\centering
\caption{Distribution}\label{tb.5.1}
\vskip .25\baselineskip
\begin{tabular}{|c||c|c|}
\hline
$t_1\backslash t_2$&$t_2^1$&$t_2^2$\\
\hline
\hline
$t_1^1$&$0.2$&$0.3$\\
\hline
$t_1^2$&$0.4$&$0.1$\\
\hline
\end{tabular}
\end{table}

Then
$$
p_1=\begin{bmatrix}
Pr(t_2^1|t_1^1),Pr(t_2^1|t_1^2)\\
Pr(t_2^2|t_1^1),Pr(t_2^2|t_1^2)\\
\end{bmatrix}
=\begin{bmatrix}
2/5,4/5\\
3/5,1/5\\
\end{bmatrix}.
$$
$$
p_2=\begin{bmatrix}
Pr(t_1^1|t_2^1),Pr(t_1^1|t_2^2)\\
Pr(t_1^2|t_2^1),Pr(t_1^2|t_2^2)\\
\end{bmatrix}
=\begin{bmatrix}
1/3,3/4\\
2/3,1/4\\
\end{bmatrix}.
$$

The payoff vectors are
$$
\begin{array}{ccl}
V_1^c&=&[a_1,a_2,a_3,a_4,b_1,b_2,b_3,b_4,c_1,c_{2},c_{3},c_{4},d_{1},d_{2},d_{3},d_{4}],\\
V_2^c&=&[e_1,e_2,e_3,e_4,f_1,f_2,f_3,f_4,g_1,g_{2},g_{3},g_{4},h_{1},h_{2},h_{3},h_{4}].
\end{array}
$$

We can calculate the payoff vectors of Action-Type Bayesian game.
$$
\begin{array}{ccl}
V_1^{AT}&=&\frac{1}{5}[2a_1+3b_1,2a_2+3b_2,2a_3+3b_3,2a_4+3b_4,4c_1+d_1,4c_2+d_2,4c_3+d_3,4c_4+d_4]\\
&:=&[\a_1,\a_2,\a_3,\a_4,\a_5,\a_6,\a_7,\a_8],\\
V_2^{AT}&=&\frac{1}{12}[4e_1+8f_1,4e_2+8f_2,4e_3+8f_3,4e_4+8f_4,9g_1+3h_1,9g_2+3h_2,9g_3+3h_3,9g_4+3h_4]\\
&:=&[\b_1,\b_2,\b_3,\b_4,\b_5,\b_6,\b_7,\b_8].\\
\end{array}
$$

According to Theorem \ref{thm5.10}
\begin{align*}
\Psi^{AT}
=\begin{bmatrix}
    -1  &   0  &   0  &   0  &   1   &  0   &  0  &   0\\
     0  &  -1  &   0  &   0  &   1   &  0   &  0  &   0\\
    -1  &   0  &   0  &   0  &   0   &  1   &  0  &   0\\
     0  &  -1  &   0  &   0  &   0   &  1   &  0  &   0\\
     0  &   0  &  -1  &   0  &   1   &  0   &  0  &   0\\
     0  &   0  &   0  &  -1  &   1   &  0   &  0  &   0\\
     0  &   0  &  -1  &   0  &   0   &  1   &  0  &   0\\
     0  &   0  &   0  &  -1  &   0   &  1   &  0  &   0\\
    -1  &   0  &   0  &   0  &   0   &  0   &  1  &   0\\
     0  &  -1  &   0  &   0  &   0   &  0   &  1  &   0\\
    -1  &   0  &   0  &   0  &   0   &  0   &  0  &   1\\
     0  &  -1  &   0  &   0  &   0   &  0   &  0  &   1\\
     0  &   0  &  -1  &   0  &   0   &  0   &  1  &   0\\
     0  &   0  &   0  &  -1  &   0   &  0   &  1  &   0\\
     0  &   0  &  -1  &   0  &   0   &  0   &  0  &   1\\
     0  &   0  &   0  &  -1  &   0   &  0   &  0  &   1\\
\end{bmatrix},~~
b^{AT}
=\begin{bmatrix}
\a_1-\b_1\\
\a_2-\b_2\\
\a_3-\b_3\\
\a_4-\b_4\\
\a_1-\b_5\\
\a_2-\b_6\\
\a_3-\b_7\\
\a_4-\b_8\\
\a_5-\b_1\\
\a_6-\b_2\\
\a_7-\b_3\\
\a_8-\b_4\\
\a_5-\b_5\\
\a_6-\b_6\\
\a_7-\b_7\\
\a_8-\b_8\\
\end{bmatrix}.
\end{align*}

Equation (\ref{eq74.1}) has a solution if and only if
$$\rank(\Psi^{AT})=\rank([\Psi^{AT},b^{AT}]),$$
which implies that
\begin{align*}
\begin{cases}
\a_1-\a_2-\a_3+\a_4-\b_1+\b_2+\b_3-\b_4=0\\
\b_2-\b_4-\b_6+\b_8=0\\
\a_1-\a_2-\a_5+\a_6=0\\
\b_1-\b_3-\b_5+\b_7=0\\
\end{cases}.
\end{align*}
Particularly,  if
$$
\begin{array}{cccccccccccl}
\a_1=\a_2,&\a_3=\a_4,&\a_5=\a_6,&\a_7=\a_8;\\
\b_1=\b_3,&\b_2=\b_4,&\b_5=\b_7,&\b_6=\b_8.\\
\end{array}
$$
The Bayesian game is an Action-Type potential game with AT potential function
$$Q^{AT}(a,t)=c^{AT}_1(a,t_1)+c^{AT}_2(a,t_2).$$
\end{exa}

\section{Dynamic Bayesian Games}

\subsection{Dynamics of Selten Bayesian Games}

For a repeated Bayesian game, its dynamics depends on two facts: (1) the strategy updating rule (SUR); (2) the conversion. Cosider a Bayesian game, a fundamental assumption in this section is, for a given conversion,  the player updates his action or type according to his SUR to optimize the corresponding conversion payoff.

\begin{dfn}
\begin{itemize}
\item[(i)]The dynamics is called asynchronous if at each time there  is only one player who is allowed to update his action or type.
\item[(ii)]The dynamics is called synchronous if at each time all players can update  actions or types synchronously.
\end{itemize}
\end{dfn}

Consider a Bayesian game $G=(N,T,A,c, Pr)$.  Assume the conversion is  S-Conversion. As the dynamics  depends on the information acquired, we assume  player $i$ knows his type $\bar{t}_i\in T_i$ and other players' actions $a_{-i}(k)$ at time $k+1$. We introduce the following asynchronous dynamics for Selten Bayesian games.
\begin{itemize}
  \item[(i)] If the SUR is asynchronous myopic best response adjustment (MBRA), then at time $k$ there is only one player (say player $i$) to update his action as follows
\begin{align} \label{6.1}
a_i(k+1)=\argmax_{a_i\in A_i}E(c_i(a_i,a_{-i}(k))|t_i=\bar{t}_i),\quad i=1,2,\cdots,n.
\end{align}
where $\bar{t}=(\bar{t}_1,\bar{t}_2,\cdots,\bar{t}_n)$ the pre-assigned type profile.
  \item[(ii)] If the SUR is  asynchronous logit response (LR),  then at time $k$ there is only one player (say player $i$) to select $a_i$ according to the following probability
\begin{align} \label{6.2}
Pr(a_i(k+1)=a_i|a(k))=\dfrac{e^{\frac{1}{T}E(c_i(a_i,a_{-i}(k))|t_i=\bar{t}_i)}}{\sum_{a'_i\in A_i}e^{\frac{1}{T}E(c_i(a'_i,a_{-i}(k))|t_i=\bar{t}_i)}},\quad i=1,2,\cdots,n.
\end{align}
\end{itemize}

We give an example to show how to get the dynamic equation for repeated Bayesian game under  S-Conversion.
\begin{exa}\label{e7.1}
Consider a Bayesian game $G=(N, A, T,c,Pr)$, where $|N|=2$, $A_1=\{a_1^1,a_1^2\}$, $A_2=\{a_2^1,a_2^2\}$,
$T_1=\{t_1^1,t_1^2\}$, $T_2=\{t_2^1,t_2^2\}$. The payoff vectors $V_1^c,V_2^c$ are
$$
\begin{array}{l}
V^c_1=[2,1,0,1,-1,1,3,-2,2,3,2,-2,3,3,-2,1],\\
V^c_2=[1,3,2,-1,2,2,1,-2,-1,0,-2,2,2,3,-1,0].\\
\end{array}
$$
The distribution shown in Table \ref{tb.6.1}.
\begin{table}[H]
\centering
\caption{Distribution}\label{tb.6.1}
\vskip .25\baselineskip
\begin{tabular}{|c||c|c|}
\hline
$t_1\backslash t_2$&$t_2^1$&$t_2^2$\\
\hline
\hline
$t_1^1$&$0.1$&$0.3$\\
\hline
$t_1^2$&$0.4$&$0.2$\\
\hline
\end{tabular}
\end{table}

Assume $t_1=t_1^1$, $t_2=t_2^2$, then we have
$$p_1({t_1^1})=(0.25,0.75)^T,\quad p_2({t_2^2})=(0.6,0.4)^T.$$
According to (\ref{4.10}), we can obtain the S-Conversion payoff vector
$$
\begin{array}{ccl}
V^S_1&=&V^c_1\d_2^1p_1({t_1^1})\\
~&=&[-0.25, 1,2.25,   -1.25],\\
V^S_2&=&V^c_2W_{[2,2]}\d_2^2p_2({t_2^2})\\
~&=&[2,2.4, 0.2, -1.2].
\end{array}
$$

Assume the SUR is  synchronous  myopic best response adjustment,
then  it is easy to get the dynamic equation as
$$
\begin{array}{l}
a_1(k+1)=\d_2[2,1,2,1]a_1(k)a_2(k),\\
a_2(k+1)=\d_2[2,2,1,1]a_1(k)a_2(k).
\end{array}
$$
Finally, we have
$$
a(k+1)=\d_4[4,2,3,1]a(k).
$$

\end{exa}

\begin{thm}
Consider a Bayesian game $G=(N, A, T,c,Pr)$. Suppose $G$ is a Selten potential game with potential function $Q^S(a)$, if each player updates its action according to asynchronous MBAR, then the dynamics converges to an S-BN-E.
\end{thm}
\noindent{\it Proof:}
Suppose the action profile at time $k$ is $a(k)$. Suppose the updating player is $i$, then

$$
\begin{array}{cccl}
&~&Q^S(a(k+1))-Q^S(a(k))~~~~~~~~~~~~~~~~~~~~~~~~~~~~~~~~~~~~~~~~~~~~~~~~~\\
&=&E(c_i(a(k+1))|t_i=\bar{t}_i)-E(c_i(a(k))|t_i=\bar{t}_i)~~~~~~~~~~~~~~~~~~~~~~~~~~~\\
&=&E(c_i(a_i(k+1),a_{-i}(k))|t_i=\bar{t}_i)-E(c_i(a_i(k),a_{-i}(k))|t_i=\bar{t}_i)\geq0.
\end{array}
$$
 Since the action profile is finite, then after finite steps, the asynchronous MBAR will converge to the potential maximizer, which is an S-BN-E of $G$.
\hfill $\Box$


\begin{thm}\label{thm7.5}
Consider a Selten Bayesian potential game $G=(N, A, T,c,Pr)$. Suppose  each player updates its action according to asynchronous LR. Then
\begin{itemize}
  \item[(i)]  The invariant distribution, denote by $\mu_T$, of asynchronous LR is
\begin{align} \label{6.1}
\mu_T(a)=\dfrac{e^{\frac{1}{T}Q^S(a)}}{\sum_{a'\in A}e^{\frac{1}{T}Q^S(a')}}.
\end{align}
  \item[(ii)] Let $\mu_0(a)=\lim_{T\rightarrow0}\mu_T(a).$ Then  the support of $\mu_0$ is equal to the set of maximizers of $Q^S(a).$
\end{itemize}
\end{thm}

\noindent{\it Proof:}
\begin{itemize}
  \item[(i)] The asynchronous LR defines an irreducible and aperiodic Markov chain on the state space $A$.
  Denote by $P(a;a')$ the probability transferring from profile $a$ to $a'$ under the dynamics of  asynchronous LR.
  To prove $\mu_T(a)$ is an invariant distribution, we only need to prove that the following detailed balance condition is satisfied
  $$\mu_T(a)P(a;a')=\mu_T(a')P(a';a),~~\forall a,a'\in A.$$
  If $a=a'$, the detailed balance condition is clearly fulfilled. If $a$ and $a'$ differ in more than one player, then $P(a;a')=P(a';a)=0.$  Hence we only consider that $a$ and $a'$ differ exactly in one player $i$. That is, $a_i\neq a'_i,~a_{-i}= a'_{-i}$. Then
  $$
\begin{array}{cccl}
  \mu_T(a)P(a;a')&=&\dfrac{e^{\frac{1}{T}Q^S(a)}}{\sum_{\bar{a}\in A}e^{\frac{1}{T}Q^S(\bar{a})}}\cdot\dfrac{e^{\frac{1}{T}c_i^S(a'_i,a_{-i})}}{\sum_{\tilde{a}_i\in A_i}e^{\frac{1}{T}c_i^S(\tilde{a}_i,a_{-i})}}\\
  &=&\dfrac{e^{\frac{1}{T}Q^S(a)}}{\sum_{\bar{a}\in A}e^{\frac{1}{T}Q^S(\bar{a})}}\cdot\dfrac{e^{\frac{1}{T}Q^S(a'_i,a_{-i})}}{\sum_{\tilde{a}_i\in A_i}e^{\frac{1}{T}Q^S(\tilde{a}_i,a_{-i})}}\\
  &=&\dfrac{e^{\frac{1}{T}Q^S(a'_i,a_{-i})}}{\sum_{\bar{a}\in A}e^{\frac{1}{T}Q^S(\bar{a})}}\cdot\dfrac{e^{\frac{1}{T}Q^S(a)}}{\sum_{\tilde{a}_i\in A_i}e^{\frac{1}{T}Q^S(\tilde{a}_i,a_{-i})}}\\
  &=&\mu_T(a')P(a';a).~~~~~~~~~~~~~~~~~~~~~~~~~~~~~~
  \end{array}
  $$
  \item[(ii)]  As $T\rightarrow 0$, the asynchronous LR converges to the myopic best reply rule. Therefore, the support of $\mu_0$ is equal to the set of maximizers of $Q^S(a).$
\end{itemize}
\hfill $\Box$

\begin{exa}\label{ebg.5.3}
Recall Example \ref{e7.1}. Assume the SUR is synchronous LR, the conversion is S-Conversion.
Assume the type assigned by nature is $t_1=t_1^1$, $t_2=t_2^2$, then the expected payoff is
$$
\begin{array}{ccccccl}
V^S_1&=&[-0.25,& 1,&2.25,& -1.25],\\
V^S_2&=&[2.000,   &2.4,& 0.20,& -1.20].
\end{array}
$$

By virtue of (\ref{6.2}), we have
\begin{align*}
Pr(a_1(k+1)=a_1|a_2(k))=
\begin{cases}
\dfrac{e^{\frac{1}{T}c_1^S(a_1=a_1^1,a_2(k)=a_2^1)}}{e^{\frac{1}{T}c_1^S(a_1=a_1^1,a_2(k)=a_2^1)}+e^{\frac{1}{T}c_1^S(a_1^2,a_2(k)=a_2^1)}}:=\a_1,\\
\dfrac{e^{\frac{1}{T}c_1^S(a_1=a_1^1,a_2(k)=a_2^2)}}{e^{\frac{1}{T}c_1^S(a_1=a_1^1,a_2(k)=a_2^2)}+e^{\frac{1}{T}c_1^S(a_1^2,a_2(k)=a_2^2)}}:=\b_1.\\
\end{cases}
\end{align*}
Similarly, we have
\begin{align*}
Pr(a_2(k+1)=a_2|a_1(k))=
\begin{cases}
\dfrac{e^{\frac{1}{T}c_1^S(a_2=a_2^1,a_1(k)=a_1^1)}}{e^{\frac{1}{T}c_1^S(a_2=a_2^1,a_1(k)=a_1^1)}+e^{\frac{1}{T}c_1^S(a_2^2,a_1(k)=a_1^1)}}:=\a_2,\\
\dfrac{e^{\frac{1}{T}c_1^S(a_2=a_2^1,a_1(k)=a_1^2)}}{e^{\frac{1}{T}c_1^S(a_2=a_2^1,a_1(k)=a_1^2)}+e^{\frac{1}{T}c_1^S(a_2^2,a_1(k)=a_1^2)}}:=\b_2.\\
\end{cases}
\end{align*}
Then
\begin{align*}
a_1(k+1)=M_1a(k),~~a_2(k+1)=M_2a(k),
\end{align*}
where $a(k)=a_1(k)a_2(k)$,
\begin{align*}
M_1=\begin{bmatrix}
\a_1&\b_1&\a_1&\b_1\\
1-\a_1&1-\b_1&1-\a_1&1-\b_1\\
\end{bmatrix}.
\end{align*}
\begin{align*}
M_2=\begin{bmatrix}
\a_2&\b_2&\a_2&\b_2\\
1-\a_2&1-\b_2&1-\a_2&1-\b_2\\
\end{bmatrix}.
\end{align*}

Finally, we have
\begin{align*}
a(k+1)=Ma(k)=(M_1*M_2)a(k),
\end{align*}
whre $*$ is the Khatra-Rao product.

Let~$T=2$, then
\begin{align*}
M_1 =
\begin{bmatrix}
    0.0067&    0.9890&    0.0067&    0.9890\\
    0.9933&    0.0110&    0.9933&    0.0110\\
\end{bmatrix}.
\end{align*}
\begin{align*}
M_2 =\begin{bmatrix}
    0.2315&    0.2315&    0.9975&    0.9975\\
    0.7685&    0.7685&    0.0025&    0.0025\\
\end{bmatrix}.
\end{align*}
\begin{align*}
M =\begin{bmatrix}
    0.0015&    0.2289&    0.0067&    0.9866\\
    0.0051&    0.7601&    0.0000&    0.0024\\
    0.2299&    0.0025&    0.9909&    0.0110\\
    0.7634&    0.0084&    0.0025&    0.0000\\
\end{bmatrix}
\end{align*}

Therefore the  steady action profile is
$$a^*(t_1^1,t_2^1)=[0.0357,0.0011,0.9336,0.0296]^T.$$
Player $1$'s steady strategy is
$$
a_1^*(t_1^1,t_2^1)=(I_2\otimes \J_2^T)a^*(t_1^1,t_2^1)=[0.0369,0.9631]^T.
$$
Player $2$'s steady strategy is
$$
a_2^*(t_1^1,t_2^1)=(\J_2^T\otimes I_2)a^*(t_1^1,t_2^1)=[0.9693,0.0307]^T.
$$

Similarly, we can conclude that
\begin{itemize}
\item If~$t_1^{\theta}(t)=t_1^1$ ~$t_2^{\theta}(t)=t_2^2$,  the  steady action profile is
$$a^*(t_1^1,t_2^2)=[0.1932,0.0680,0.5466,0.1922]^T.$$
Player $1$'s steady strategy is
$$
a_1^*(t_1^1,t_2^2)=(I_2\otimes \J_2^T)a^*(t_1^1,t_2^2)=[0.2612,0.7388]^T.
$$
Player $2$'s steady strategy is
$$
a_2^*(t_1^1,t_2^2)=(\J_2^T\otimes I_2)a^*(t_1^1,t_2^2)=[0.7398,0.2602]^T.
$$

\item
If~$t_1^{\theta}(t)=t_1^2$ ~$t_2^{\theta}(t)=t_2^1$,  the  steady action profile is
$$a^*(t_1^2,t_2^1)=[ 0.2359,0.7556,0.0020,0.0064]^T.$$
Player $1$'s steady strategy is
$$
a_1^*(t_1^2,t_2^1)=(I_2\otimes \J_2^T)a^*(t_1^2,t_2^1)=[ 0.9915,0.0085]^T.
$$
Player $2$'s steady strategy is
$$
a_2^*(t_1^2,t_2^1)=(\J_2^T\otimes I_2)a^*(t_1^2,t_2^1)=[  0.2379,0.7621]^T.
$$

\item
If~$t_1^{\theta}(t)=t_1^2$ ~$t_2^{\theta}(t)=t_2^2$,  the  steady action profile is
$$a^*(t_1^2,t_2^2)=[ 0.0858,0.9109,0.0003,0.0030]^T.$$
Player $1$'s steady strategy is
$$
a_1^*(t_1^2,t_2^2)=(I_2\otimes \J_2^T)a^*(t_1^2,t_2^2)=[  0.9967,0.0033]^T.
$$
Player $2$'s steady strategy is
$$
a_2^*(t_1^2,t_2^2)=(\J_2^T\otimes I_2)a^*(t_1^2,t_2^2)=[ 0.0861,0.9139]^T.
$$
\end{itemize}

\end{exa}
%
%
%
%
%
%
%

\subsection{Dynamics of Action-Type Bayesian Games}
Consider a Bayesian game $G=(N,T,A,c, Pr)$. Assume the conversion is  AT-Conversion. Unlike  Selten Bayesian game, the SUR of AT Bayesian game consists of action updating rule (AUR) and type updating rule (TUR). If the AUR of player $i$ is $f_i$, and the TUR of player $i$ is $g_i$, then the SUR is
\begin{align}
\begin{array}{l}
\begin{cases}
a_i(k+1)=f_i(a_1(k),\cdots,a_n(k),t(k)),\\
t_i(k+1)=g_i(a_1(k),\cdots,a_n(k),t(k)).
\end{cases}
\end{array}
\end{align}
The player can update its action and type, concurrently or separately. Assume in the decision-making process player $i$ knows his type $t_i(k)$ and other player's actions $a_{-i}(k)$ at time $k+1$, but he doesn't know  other player's type $t_{-i}(k)$.

Firstly, we design different asynchronous MBARs for Action-Type Bayesian games, which are shown as follows:
\begin{itemize}
  \item[(i)] The SUR is called asynchronous concurrent MBRA (C-MBRA), if at time $k$ the updating player $i$ chooses his action and type concurrently, which is shown as follows
\begin{align} \label{6.1}
a_i(k+1)\ltimes t_i(k+1)=\argmax_{(a_i,t_i)\in A_i\times T_i}c_i^{AT}(a_i,t_i,a_{-i}(k)),\quad i=1,2,\cdots,n.
\end{align}
 \item[(ii)] The SUR is called asynchronous separate MBRA (S-MBRA), if at time $k$ the updating player $i$ either chooses his action  as follows
\begin{align} \label{6.1}
a_i(k+1)=\argmax_{a_i\in A_i}c_i^{AT}(a_i,t_i(k),a_{-i}(k)),\quad i=1,2,\cdots,n.
\end{align}
or selects type as follows
\begin{align} \label{6.1}
t_i(k+1)=\argmax_{t_i\in T_i}c_i^{AT}(t_i,a_i(k),a_{-i}(k)),\quad i=1,2,\cdots,n.
\end{align}
\end{itemize}

\begin{thm}\label{thm6.4}
Consider an Action-Type Bayesian potential game $G=(N, A, T,c,Pr)$. If each player updates its action according to asynchronous C-MBAR or S-MBAR, then the dynamics converges to AT-BN-E.
\end{thm}
\noindent{\it Proof:}
Suppose the action-type profile at time $k$ is $(a(k),t(k))$. Suppose the updating player is $i$ at time $k$ using asynchronous C-MBAR, then
$$
\begin{array}{cccl}
&~&Q^{AT}(a(k+1),t(k+1))-Q^{AT}(a(k),t(k))~~~~~~~~~~~~~~~~~~~~~~~~~~~~~~\\
&=&c^{AT}_i(a(k+1),t_i(k+1))-c^{AT}_i(a(k),t_i(k))~~~~~~~~~~~~~~~~~~~~~~~~~~~~\\
&=&c^{AT}_i(a_i(k+1),a_{-i}(k),t_i(k+1))-c^{AT}_i(a_i(k),a_{-i}(k),t_i(k))\geq0,\\
\end{array}
$$
where $Q^{AT}$ is the potential function of the Action-Type Bayesian potential game $G$.
 Since the action profile is finite,  after finite steps, the asynchronous MBAR will converge to the potential maximizer, which is an AT-BN-E of $G$. The proof to asynchronous S-MBAR is the same, so we omit the details.

\hfill $\Box$

Similarly, we can design asynchronous concurrent LR (C-LR) and separate LR (S-LR) for Action-Type Bayesian games, which are shown as follows:
\begin{itemize}
  \item[(i)] The SUR is called asynchronous C-LR, if at time $k$ the updating player $i$ chooses his action and type concurrently according to the following probability
\begin{align} \label{6.1}
Pr\big(a_i(k+1)=a_i,t_i(k+1)=t_i|a(k),t_i(k)\big)=\dfrac{e^{\frac{1}{T}c_i^{AT}(a_i,t_i,a_{-i}(k))}}{\sum_{(a'_i,t'_i)\in A_i\times T_i}e^{\frac{1}{T}c_i^{AT}(a'_i,t'_i,a_{-i}(k))}},\quad i=1,2,\cdots,n.
\end{align}
\item[(ii)] The SUR is called asynchronous S-LR, if at time $k$ the updating player $i$ either chooses his action according to the following probability
\begin{align} \label{6.1}
Pr\big(a_i(k+1)=a_i|a(k),t_i(k)\big)=\dfrac{e^{\frac{1}{T}c_i^{AT}(a_i,t_i(k),a_{-i}(k))}}{\sum_{a'_i\in A_i}e^{\frac{1}{T}c_i^{AT}(a'_i,t_i(k),a_{-i}(k))}},\quad i=1,2,\cdots,n.
\end{align}
or chooses his type according to the following probability
\begin{align} \label{6.1}
Pr\big(t_i(k+1)=t_i|a(k),t_i(k)\big)=\dfrac{e^{\frac{1}{T}c_i^{AT}(t_i,a_i(k), a_{-i}(k))}}{\sum_{t'_i\in T_i}e^{\frac{1}{T}c_i^{AT}(t'_i,a_i(k),a_{-i}(k))}},\quad i=1,2,\cdots,n.
\end{align}
\end{itemize}

\begin{thm}\label{thm6.5}
Consider an  Action-Type Bayesian potential game $G=(N, A, T,c,Pr)$. Suppose  each player updates his action according to asynchronous C-LR or S-LR. Then
\begin{itemize}
  \item[(i)]  The invariant distribution, denote by $\nu_T$, of asynchronous LR (either C-LR or S-LR) is
\begin{align} \label{6.1}
\nu_T(a,t)=\dfrac{e^{\frac{1}{T}Q^{AT}(a,t)}}{\sum_{(a',t')\in A\times T}e^{\frac{1}{T}Q^{AT}(a',t')}}.
\end{align}
  \item[(ii)] Let $\nu_0(a,t)=\lim_{T\rightarrow0}\nu_T(a,t).$ Then  the support of $\nu_0$ is equal to the set of maximizers of potential function $Q^{AT}(a,t).$
\end{itemize}
\end{thm}
\noindent{\it Proof:} We only prove (i), because the proof of (ii) is similar to the proof of Theorem \ref{thm7.5}.

The asynchronous C-LR defines an irreducible and aperiodic Markov chain on the state space $A\times T$.
Denote by $P(a,t;a',t')$ the probability transferring from profile $(a,t)$ to $(a',t')$ under the dynamics of  asynchronous C-LR.
To prove $\nu_T(a,t)$ is an invariant distribution, we only need to prove that the following detailed balance condition is satisfied
$$\nu_T(a,t)P(a,t;a',t')=\nu_T(a')P(a',t';a,t),~~\forall a,a'\in A,~t,t'\in T.$$

If $a=a',t=t'$, the detailed balance condition is clearly fulfilled. If $a,a'$ or $t,t'$ differ in more than one player, then $P(a,t;a',t')=P(a',t';a,t)=0.$  Hence we only consider that $a,a'$ or $t,t'$ differ exactly in one player $i$. That is, $a_i\neq a'_i,~a_{-i}= a'_{-i};~t_i\neq t'_i,~t_{-i}= t'_{-i}$. Then,
$$
\begin{array}{cccl}
  \nu_T(a,t)P(a,t;a',t')&=&\dfrac{e^{\frac{1}{T}Q^{AT}(a,t)}}{\sum_{(\bar{a},\bar{t})\in A\times T}e^{\frac{1}{T}Q^{AT}(\bar{a},\bar{t})}}\cdot\dfrac{e^{\frac{1}{T}c_i^{AT}(a'_i,a_{-i},t'_i)}}{\sum_{(\tilde{a}_i,\tilde{t}_i)\in A_i\times T_i}e^{\frac{1}{T}c_i^{AT}(\tilde{a}_i,\tilde{t}_i,a_{-i})}}\\
  &=&\dfrac{e^{\frac{1}{T}Q^{AT}(a,t)}}{\sum_{(\bar{a},\bar{t})\in A\times T}e^{\frac{1}{T}Q^{AT}(\bar{a},\bar{t})}}\cdot\dfrac{e^{\frac{1}{T}Q^{AT}(a'_i,a_{-i},t'_i,t_{-i})}}{\sum_{(\tilde{a}_i,\tilde{t}_i)\in A_i\times T_i}e^{\frac{1}{T}Q^{AT}(\tilde{a}_i,\tilde{t}_i,a_{-i},t_{-i})}}\\
  &=&\dfrac{e^{\frac{1}{T}Q^{AT}(a'_i,a_{-i},t'_i,t_{-i})}}{\sum_{(\bar{a},\bar{t})\in A\times T}e^{\frac{1}{T}Q^{AT}(\bar{a},\bar{t})}}\cdot\dfrac{e^{\frac{1}{T}Q^{AT}(a,t)}}{\sum_{(\tilde{a}_i,\tilde{t}_i)\in A_i\times T_i}e^{\frac{1}{T}Q^{AT}(\tilde{a}_i,\tilde{t}_i,a_{-i},t_{-i})}}\\
  &=& \nu_T(a',t')P(a',t';a,t).~~~~~~~~~~~~~~~~~~~~~~~~~~~~~~~~~~~~~~~~~~~~~~
\end{array}
$$
The proof for asynchronous S-LR is the same, so it is omitted.

\hfill $\Box$

We give an example to show how to get the dynamic equation for repeated Bayesian game under  AT-Conversion.
\begin{exa}
Recall the Bayesian  game $G=(N, A, T,c,Pr)$ in Example \ref{e5.11}. Suppose
$$
\begin{array}{cccccccccccl}
\a_1=-1,&\a_2=3,&\a_3=-4,&\a_4=2,&\a_5=-1,&\a_6=5,&\a_7=1,&\a_8=2;\\
\b_1=7,&\b_2=5,&\b_3=6,&\b_4=4,&\b_5=-2,&\b_6=3,&\b_7=1,&\b_8=2.\\
\end{array}
$$
According to Example \ref{e5.11}, $G$ is   Action-Type  potential.

Assume the SUR is  asynchronous C-MBRA, player $1$ updates his action and type first and then player $2$ updates his action and type. It is easy to get the dynamic equation as
$$
\begin{array}{l}
\begin{cases}
t_1(2k+1)a_1(2k+1)=\d_4[4,3,4,3,4,3,4,3]t_1(2k)a(2k),\\
t_2(2k+1)a_2(2k+1)=t_2(2k)a(2k),~~~k=1,2,\cdots,
\end{cases}
\end{array}
$$
and
$$
\begin{array}{l}
\begin{cases}
t_2(2k)a_2(2k)=\d_4[1,1,1,1,1,1,1,1]t_2(2k-1)a(2k-1),\\
t_1(2k)a_1(2k)=t_1(2k-1)a(2k-1),~~~k=1,2,\cdots.
\end{cases}
\end{array}
$$

Denote by $x(k):=t(k)\ltimes a(k)$. It is obvious that the fixed points of asynchronous C-MBRA dynamics are
$$x_1^f=(t_1^2,t_2^1,a_1^1,a_2^1),~~~x_2^f=(t_1^2,t_2^1,a_1^2,a_2^1).$$
We can verify that both $x_1^f$ and $x_2^f$ are the AT-BN-E of $G$.


\end{exa}

\section{Conclusion}

The main purpose of this paper is to provide a fundamental framework for finite Bayesian games, using STP expression of finite games.

Three conversions from incomplete information to complete information are discussed, which are Harsanyi, Selten, and AT conversions respectively. Formulas are obtained for three conventions respectively.
In addition,  Bayesian potential game is also investigated.
Finally, the evolutive equations for dynamic Bayesian games are also provided according to SURs and conversions.

There are many problems remain for further investigation. For instance, the following are some challenging topics:

\begin{itemize}
\item[(i)] A generalization: State-based game. (As time-varying types.)
\item[(ii)] Deal with all history based knowledge SUR. (Truncation)
\item[(iii)] Learning-based optimization.
\item[(iv)] Calculating Bayesian Nash equilibrium.
\item[(v)] Dynamic Bayesian Nash equilibrium.
\item[(vi)] Application of Potential Bayesian game.
\end{itemize}

\end{document}